\documentclass[12pt]{article}
\usepackage{amsmath,amsfonts,amsthm,latexsym}

  \def\CC{{\mathbb C}} \def\DD{{\mathbb D}}
 \def\FF{{\mathbb F}}  \def\HH{{\mathbb H}}
   
 \def\NN{{\mathbb N}}  
 \def\RR{{\mathbb R}}  
     \def\ZZ{{\mathbb Z}}

\def\Si{\Sigma}

\def\De{\Delta}

\def\Ga{\Gamma}

\def\cA{{\cal A}}

\def\cE{{\cal E}}    
\def\cF{{\cal F}}  \def\cL{{\cal L}} \def\cR{{\cal R}}  \def\cY{{\cal Y}}

\newtheorem{theo}{Theorem}

\newtheorem{corol}[theo]{Corollary}
\newtheorem{claim}{Claim}

\newtheorem{lemm}{Lemma}[section]
\newtheorem{coro}[lemm]{Corollary}
\newtheorem{defi}[lemm]{Definition}
\newtheorem{prop}[lemm]{Proposition}

\newtheorem{rema}[lemm]{Remark}
\newenvironment{demo}{{\bf Proof: }}{\hfill$\Box$\medskip}

\title{\bf Statistical Behaviour of the Leaves of Riccati Foliations
\thanks{Partially supported by CNRS, France, and CONACYT 61317, 28541-E
and 28491-E, Mexico.}}
\author{Ch. Bonatti, X. G\'omez-Mont and R. Vila-Freyer }

\date{\today}

\begin{document}
\maketitle

\begin{abstract} We introduce the geodesic
flow on the leaves of a holomorphic foliation with leaves of
 dimension 1 and hyperbolic, corresponding to the
 unique complete metric of curvature -1 compatible
 with its conformal structure.
 We do these for the foliations associated to Riccati
 equations, which
 are
the projectivisation of the solutions of a
linear ordinary differential equations over a finite
Riemann surface of hyperbolic type $S$,
and may be described by a representation
$\rho:\pi_1(S) \rightarrow GL(n,\CC)$.
We give conditions under which the foliated geodesic flow
has a generic repellor-attractor statistical dynamics.
That is, there are measures $\mu^-$ and $\mu^+$
such that for almost any initial condition with respect
to the Lebesgue measure class the statistical average of
the foliated geodesic flow converges
for negative time to $\mu^-$ and for positive time to $\mu^+$ 
(i.e. $\mu^+$ is the unique SRB-measure and its basin has total Lebesgue measure).
These measures are ergodic
with respect to the foliated geodesic flow. These measures
are also invariant under a foliated horocycle flow and 
they project to a harmonic measure for the Riccati foliation,
which plays the role of an attractor for the statistical behaviour of
the leaves of the foliation.
\end{abstract}

\section*{Introduction}
\vskip 5mm

The objective of this work is to propose a method
for
 understanding the statistical properties
of the leaves of a holomorphic foliation,
and which we carry out for
a simple class of holomorphic foliations:
those obtained from the solutions of Riccati
Equations.
The method consists in using
the canonical
metric of curvature -1 that the
leaves have as Riemann surfaces,
the Poincar\'e metric, and then to flow
along foliated geodesics.
One is  interested in understanding the statistics
of this foliated geodesic flow.
In particular, in determining if
the foliated geodesic flow has an SRB-measure (for Sina\"\i~, Ruelle and
Bowen [21], [20], [5]),
or physical measure,
which means that
a set of geodesics of positive Lebesgue measure
have a convergent time statistics,
which is shared by all the geodesics
in this set, called the basin of
attraction of the SRB-measure.
The SRB-measure is the spatial measure
describing this common time statistics
of a significant set of geodesics.
One then shows that the SRB-measure
is invariant also under a foliated horocycle flow ([2]) and
 the projection of the SRB-measure to the
$\CC P^{n-1}$-bundle over $S$  
is a harmonic measure for the Riccati foliation;
in fact, the harmonic measures are in 1-1 correspondance 
with the measures simultaneously invariant
by the foliated geodesic and a horocycle flow
([1], [17], [18]]).

\vskip 2mm
The approach of using harmonic measures to understand the statistical
behaviour of the leaves of a foliation 
started with the work of Garnett([11]) who proved existence 
of harmonic measures for
regular foliations in compact manifolds, containing 
statistical properties of the behaviour of the leaves of the foliation. 
In this work we are dealing
with singular foliations in compact manifolds 
(obtained by compactifying 
the Riccati foliation
with a linear model with singularities over each puncture of $S$),
 which introduces the difficulty 
that the support of the measures could
be contained in the singular set. 
Our conclusions are related to Fornaess and Sibony's harmonic currents in $\CC P^2$
([8], [9] and [10]),
where they show existence and uniqueness 
of harmonic currents
using $\bar \partial$-methods for the
generic foliations in $\CC P^2$. 
Their result does not include Riccati foliations in $\CC P^2$,
 since these
have some tangent lines (corresponding to the punctures of $S$) and 
a non-hyperbolic singular point (arising from the blow down to $\CC P^2$).
 Our work is also related to Deroin and Kleptsyn [7],
where they use foliated Brownian motion and heat flow
instead of the foliated
geodesic and horocycle flows for
non-singular transversely holomorphic foliations 
in compact manifolds
to obtain a finite number of
attracting harmonic measures and the negativity of
the Lyapunov exponent.

\vskip 2mm
The Riccati equations
are projectivisations of linear ordinary
differential equations over
a finite hyperbolic Riemann surface $S$
(i.e. compact minus a finite number of points
and with universal cover the upper half plane).
Locally they have the form

$$ \frac{dw}{  dz} = A(z)w\hskip 5mm,\hskip5mm w\in
\CC^n \ , \ z \in \CC \hskip 5mm, \hskip 5mm A:\CC \to Mat_{n,n}(\CC)
$$
with $A$ holomorphic.
These equations may be equivalently defined by giving the
monodromy representations

$$
\tilde \rho : \pi_1(S,z_0) \rightarrow GL(n,\CC)
\hskip 1cm , \hskip 1cm
\rho : \pi_1(S,z_0) \rightarrow PGL(n,\CC)
\eqno(1)$$
and suspending them, to obtain flat $\CC^n$ and $\CC P^{n-1}$
bundles over $S$

$$E_{\tilde \rho} \rightarrow S, \hskip 1cm,\hskip1cm
  M_{ \rho} \rightarrow S.
  \eqno (2)$$
The graphs of the local flat sections of these bundles
are the `solutions' to the linear differential equation
defined by the monodromy $(1)$  and
define holomorphic foliations $
\cF_{\tilde \rho}$ and
$\cF_{ \rho}$
of $
E_{\tilde \rho}
$ and $  M_{ \rho}$
whose leaves $\cL$
project as a covering
to the base surface $S$.

\vskip 2mm
Introduce to the finite hyperbolic Riemann surface $S$
the Poincar\'e metric, to the
unit tangent bundle $q:T^1S\rightarrow S$
the geodesic flow $\varphi:T^1S\times \RR \rightarrow T^1S$
and the Liouville measure   $dLiouv$
 (hyperbolic area element in $S$
and Haar measure on $T^1_pS$, normalised to volume 1).
 We may introduce on the leaves $\cL$ of the
 foliations $\cF_{\tilde \rho}$ and
 $\cF_\rho$ the Poincar\'e metric,
 which is the pull back of
 the Poincar\'e metric of $S$ by the covering map
 $q:\cL \rightarrow S$.
 The unit tangent bundle
 $T^1_{{\cF}_{\tilde\rho}}$
 to the
 foliation $\cF_{\tilde\rho}$
 in $E_{\tilde \rho}$
is canonically isomorphic
to
 the vector bundle
$q^* E_{\tilde\rho}$
over $ T^1S$,
that we denote by  $
E $. In the same way the unit tangent bundle $T^1_{\cF_\rho}$ of the foliation $\cF_\rho$ is canonically identified to the projectivisation $Proj(E)$ of the vector bundle $E$ over $T^1S$. 
Introduce
on $E$ and on $Proj(E)$
the foliated
geodesic flows $\tilde \Phi$ and $\Phi$
(see $(2.2)$), obtained by
flowing along the foliated geodesics.
Introduce also on $E$ a continuous Hermitian inner product $|\;.\;|_v$.

\vskip 2mm
Given a vector $v \in T^1S$ we have
the geodesic
$$\RR \rightarrow T^1S
\hskip 1cm,\hskip 1cm
t \rightarrow \varphi(v,t)$$ determined by
the initial condition $v$
and given  $w_0 \in E_v$ we
 also have the foliated geodesic
$$\RR \rightarrow E
\hskip 1cm,\hskip 1cm
t \rightarrow \tilde \Phi(w_0,t)$$
which is the solution to
the linear differential equation
defined by $(1)$ along the foliated geodesic
determined by $v$ and $w_0$. The function
$$t \rightarrow | \tilde \Phi (w_0,t) |_{\varphi(v,t)}$$
describes the type of growth of the solution of $(1)$
along the geodesic $\gamma_v$ with initial
condition $w_0 \in E_v$
and  the function
$$t \rightarrow \frac{| \tilde \Phi (w_1,t) |_{\varphi(v,t)}}{
 | \tilde \Phi (w_2,t) |_{\varphi(v,t)}}$$
describes the relative growth of the solution of $(1)$
along the geodesic $\gamma_v$ with initial
condition $w_1 \in E_v$
with respect to the growth of the solution of $(1)$
along $\gamma_v$ with the initial condition $w_2 \in E_v$.

\vskip 2mm
We say that the Riccati equation has a
{\bf  section of largest expansion} $\sigma^+$
if
for Liouville almost any point $v$ on $T^1S$ we may
measurably define
a splitting $E_v = F_v \oplus G_v$ by linear spaces,
which is invariant by the foliated geodesic flow
$\tilde \Phi$
with $F_v$ of dimension 1 and
with the property that the map
$ t \rightarrow
\tilde \Phi(w_1,t)$
with initial condition
$w_1 \in F_v$ grows more rapidly
than the maps
$ t \rightarrow
\tilde \Phi(w_2,t)$
for any $w_2 \in G_v$. That is,
 for almost any $v\in T^1S$, for any compact set $K\subset T^1S$
and for any sequence $(t_n)_{n\in\NN}$ of times
such that $\varphi(v,t_n)\in K$ and
$\lim_{n\to\infty}t_n=+\infty$, one has:
$$\lim_{n\to \infty} \frac{ |
\tilde \Phi(w_1,t_n)|_{\varphi(v,t_n)}}{ |
\tilde \Phi(w_2,t_n)|_{\varphi(v,t_n)}} = \infty,
\quad\mbox{for all non-zero}\quad w_1\in F_v,
\quad\mbox{and}\quad w_2\in G_v.
$$
So the section of largest expansion is defined as
  $\sigma^+:= Proj(F):T^1S \rightarrow Proj(E)$.
Similarly, we may define  a
{\bf   section $\sigma^-$ of largest contraction} (see $(3.1)$).

  \vskip 2mm
An elementary argument of Linear Algebra
suggests that a section  $\sigma^+=Proj(F)$
of largest expansion
 is attracting all the points in $Proj(E) -Proj(G)$
as they flow according to 
the action of the foliated
 geodesic flow $ \Phi$.
In fact, we prove:

\vskip 2mm
\begin{theo}
Let $S$ be a finite hyperbolic  Riemann surface
and $\tilde \rho:\pi_1(S,z_0) \rightarrow  GL(n,\CC)$ a representation
having a section
 $\sigma^+$
of largest expansion,
then $\mu^+=\sigma^+_*(dLiouv)$ is a $\Phi-$invariant
ergodic measure on $T^1\cF_\rho$
which is an SRB-measure 
for the foliated geodesic flow $\Phi$
of the Riccati equation, whose basin has total Lebesgue
measure in $T^1\cF_\rho$.
Similarly,
if $\sigma^-$ is the
section of largest contraction,
then
 $\mu^-=\sigma^-_*(dLiouv)$ is a $\Phi-$invariant
ergodic measure which is an SRB-measure whose basin has total Lebesgue
measure in $T^1\cF_\rho$, for negative times.
\label{theo1}
\end{theo}

\vskip 2mm
In the case that both $\sigma^\pm$ exist,
the foliated geodesic flow has a very simple `north to
south pole dynamics': almost everybody is born in
$\mu^-$ and is dying on $\mu^+$.
If the sections $\sigma^\pm$ are continuous
disjoint sections
defined on all $T^1S$
then it is easy to imagine this north to south pole dynamics
(see section 7 for an example).
If $\sigma^\pm$ are only measurable, then 
they describe more subtle phenomena.

\vskip 2mm
The Lyapunov exponents measure
the exponential rate of growth (for the metric $|\;.\;|_v$ in the vectorial fibers) of the solutions
of the linear equation along the geodesics
(definition 4.2):
$$\lim_{t \rightarrow \pm \infty}\frac{ 1 }{ t} \log |\tilde \Phi(w_0,t)|_{\varphi(v,t)}.$$

\vskip 2mm
Let $S$
be a finite hyperbolic  Riemann surface,
$\tilde\rho\colon\pi_1(S,z_0)\to GL(n,\CC)$  a
representation and $E$ the previously constructed bundle.
The association of initial conditions to final conditions
for the linear equation in $E$ 
over the geodesic flow of $S$, after a measurable trivialisation of the bundle,
gives rise to a 
measurable multiplicative cocycle
over the geodesic flow on $T^1S$
$$\tilde A:T^1S \times \RR \longrightarrow GL(n,\CC)$$
(see $(2.4))$.
The integrability
condition
$$
\int_{T^1{S}}log^+\|  \tilde A_{\pm 1} \| dLiouv < +\infty,
\eqno(3)$$
where $\|\ \|$ is the operator norm and $\tilde A_t :=\tilde
A(\cdot,t)$, 
asserts that the amount of expansion of $\tilde A_\pm 1$ is Liouville integrable.

\vskip 2mm
As a consequence of  the multiplicative Ergodic Theorem of Oseledec
applied to the foliated geodesic flow
we obtain:

\vskip 2mm
\begin{corol} Let $S$
be a finite hyperbolic  Riemann surface,
$\tilde\rho\colon\pi_1(S,z_0)\to GL(n,\CC)$  a
representation and let
$\tilde A$
be the measurable multiplicative cocycle
over the geodesic flow on $T^1S$
satisfying the integrability
condition $(3)$, then:
\begin{itemize}
\item
 The Lyapunov exponents
 $\lambda_1<\cdots<\lambda_k$
 of $\tilde \Phi$ are well
 defined and are constant on a subset of $T^1S$
 of total Liouville measure. Denote by 
 $F_i(v)$ the corresponding Lyapunov spaces.
\item For every $i\in\{1,\dots,k\}$, $\lambda_{k+1-i}=-\lambda_i$ and $dim(F_{k+1-i})=dim(F_i)$.
\item If $dim F_k=1$, denote by $\sigma^+$
the section corresponding to $F_k$ and
$\sigma^-$ the section correponding to $F_1$,
then $\sigma^\pm$ are sections of largest expansion
and contraction, respectively.
\end{itemize}
\end{corol}

\vskip 2mm
From now on by  the {\em Lyapunov exponents
of the linear equation obtained from the representation $\tilde \rho$} we will
understand the Lyapunov exponents of the above multiplicative
cocycle $\tilde A$ over the geodesic flow on $T^1S$ obtained
from the foliated geodesic flow on $E$
and satisfying the integrability condition $(3)$.
The relationship between the section of largest expansion
and the Lyapunov exponents
is:

\vskip 2mm
\begin{theo}
 Let $S$
be a finite hyperbolic  Riemann surface,
$\tilde\rho\colon\pi_1(S,z_0)\to GL(n,\CC)$  a
representation satisfying the integrability condition $(3)$,
then
there exists a section of largest expansion if and only if
the largest Lyapunov exponent is positive and simple,
if and only if the smallest Lyapunov exponent is negative and simple,
and 
if and only if there is a section of largest contraction.
\label{theo2}
\end{theo}

So a section of largest expansion
is an extension for  non-integrable cocycles $\tilde A$
of the notion of having a simple largest Lyapunov exponent.
We  give an example of this in section 6.

\vskip 2mm
In order to apply Oseledec's   Theorem,
the prevailing hypothesis is the integrability
condition $(3)$.
 This condition is always satisfied
 if the base Riemann surface is compact, and
more generally:

\vskip 2mm
\begin{theo}
If $S$ is a  finite hyperbolic  Riemann surface,
$\tilde \rho$ a representation $(1)$ then
the multiplicative cocycle $\tilde A$ satisfies the integrability
condition $(3)$ if and only if
the monodromy $\tilde \rho$ around each of the punctures of $S$
corresponds to a matrix with all its eigenvalues
of norm 1.
\label{t.6}
\end{theo}

\vskip 2mm
We then develop two kinds of examples: The ping-pong  or Schottky
monodromy
representations in $SL(2,\CC)$ and the
canonical representation obtained from
the representation 
$$\rho_{can}:\pi_1(S,z_0)
 \rightarrow SL(2,\RR) \subset SL(2,\CC)
\eqno (4)$$
on the universal covering of the
surface. We obtain:

\begin{theo} Let $S$ be a finite
hyperbolic   non-compact Riemann surface and $\rho\colon\pi_1(S
,z_0)\to GL(2,\CC)$ an injective
representation onto a
Schottky group, then
there are sections $s^+$ and $s^-$ of largest expansion
and contraction defined and continuous
on a   subset of $T^1S$
of full Liouville  measure.
\label{t.Schottky}
\end{theo}

\vskip 2mm
\noindent
It follows from Theorem 4 that the Schottky representations
in Theorem 5
do not satisfy
the integrability condition
$(3)$,
but
we obtain that
there are still sections of largest
expansion and contraction.
We think that
the Lyapunov exponents are in this case 
 $\pm \infty$. In fact, we can prove this assertion for
specific Schottky representations.

\begin{theo} For any finite hyperbolic Riemann surface $S$
 the foliated geodesic flow associated to  the canonical
representation $(4)$
 admits sections
of largest expansion and contraction defined
and smooth on all $T^1S$.
Moreover, for Lebesgue almost any point of $Proj(E)$
the foliated geodesic starting at this point has
$\mu^+$ as its positive statistics and
$\mu^-$ as its negative statistics
(that is, $\mu^+$ is the unique SRB-measure and its
basin has total Lebesgue measure, and similarly
$\mu^-$ for negative time).

If $S$ is compact then $\sigma^+(T^1S)$
 is a hyperbolic attractor
 and  $\sigma^-(T^1S)$
 is a hyperbolic repellor with
 basins of attraction
$T^1\cF - \sigma^\mp(T^1S)$.
\label{t.tautologic}
\end{theo}

\vskip 2mm
The  statements and arguments presented here   
extend 
to the case when the representation 
 $\rho:\pi_1(S) \rightarrow PGL(n,\CC)$  does not admit 
a lifting to a representation in $GL(n,\CC)$.

\vskip 2mm
Restricting now to $n=2$ or $3$,  assuming the integrability condition $(3)$ and
that the representation $\rho$ does not leave invariant any probability measure
(which is a generic condition on $\rho$), it follows from Theorem 3 in [2]
that the SRB-measure of the geodesic flow $\mu^+$ is the unique measure invariant
under the foliated stable horocycle flow $H^{uu}_\rho$ that projects to the Liouville
measure on $S$. Furthermore, it follows from the arguments in [1] and [17] that 
the projection to $M_\rho$ of $\mu^+$ is the unique harmonic measure $\nu$ of the
Riccati foliation  ${\cF}_{\rho}$ that projects to the Liouville
measure on $S$. It is shown in [2]
 that  $\nu$
describes effectively the statistical behaviour
of the leaves of the foliation $\cF_\rho$:
For any compact set $K \subset M_\rho$, for any
sequence $(x_n \in K)_{n\in \NN}$ and any sequence
of real numbers $(r_n)_{n\in \NN}$ tending to $+\infty$ the
family of probability measures 
$\nu_{r_n}(x_n)$ obtained by normalizing
the area element on the disk
$D_{r_n}(x_n)$ in the leafwise Poincar\'e metric
converges towards $\nu$ for the weak topology when $n$
tends to $+\infty$.
If $S$ is compact, then the integrability condition $(3)$ is always satisfied 
and the condition of projecting to the Liouville measure on $S$ is satisfied
automatically by Hedlund's Theorem [13].

 \vskip 2mm
If $S$ is a compact hyperbolic Riemann surface,
then the foliated geodesic flow is a linear or projective
multiplicative cocycle over a hyperbolic dynamical system.
This led us to think that it could be possible to adapt
Fustenberg's theory of the existence of a positive
Lyapunov exponent for random products of matrices.
This has been carried out in [3]. 
It seems possible that using a generalization of [3] found in [4] 
(simplicity of the Lyapunov spectrum) and [22] (generalization
for linear cocycles over non-uniform hyperbolic measures),
one may extend the above mentioned results for $n\geq 4$ and $S$ a finite hyperbolic
Riemann surface.

\vskip 2mm
This paper is organised as follows. In section 1
we recall the Riccati equations and in section 2 we set up the
foliated geodesic flow on Riccati equations.
In section 3 we introduce SRB-measures and prove Theorem 1.
In section 4 we prove Corollary
 2  and Theorem 3.
In section 5 we
prove Theorem 4. In sections 6 and 7
we describe the examples, proving Theorems 5 and 6.

 \vskip 2mm
\section{The Riccati Equation\label{s.Riccati}}

\subsection{Linear Ordinary Differential Equations\label{ss.edol}}

{ The classical  linear
ordinary differential equation} is

$$\frac{dw}{  dz}  = A(z)w \hskip 1cm, \hskip 1cm z\in \CC, \ w\in \CC^n
\eqno (1.1)$$
where $A(z)$ is a matrix of rational functions
(see [6]).
The
fundamental property of this equation is that locally in
$z$ we can find a basis of independent solutions of (1.1)
which accept analytic continuation to the universal
covering space of $S := \bar\CC - \hbox{poles}(A)$ as
holomorphic vector valued functions $w$
satisfying the monodromy relation:

$$ w(T_\gamma(z)) = \tilde \rho (\gamma)( w(z)) \hskip 1cm
,\hskip 1cm \gamma \in \pi_1(S,z_0)$$
where $T_\gamma$ is the covering transformation
corresponding to the close loop $\gamma$
and
$$\tilde \rho :\pi_1(S,z_0) \rightarrow GL(n,\CC)
\eqno (1.2)$$
is the monodromy representation of the equation.
The linear automorphism
$\tilde \rho (\gamma):\CC^n \rightarrow \CC^n$
contains the information of how the initial conditions
are transformed to final conditions by solving $(1.1)$
along the closed loop $\gamma$ based at $z_0$.

\vskip 2mm
Another classical construction  of  linear ordinary differential
equations is the suspension ([16]).
Assume given a hyperbolic
Riemann surface $S$
and a representation $(1.2)$.
We  construct from these data a vector bundle $E_{\tilde\rho}$ over
$S$ and an equation of type $(1.1)$.
Let $\HH^+$ be the upper half plane, considered as the universal
covering space of $S$, with covering transformations
$(4)$
giving rise to the canonical representation 
$\tilde \rho _{can}$ of the
fundamental group of $S$.
Consider the trivial bundle $\tilde E :={\HH^+\times \CC^n}$
on the upper half plane $\HH^+$
and the $\pi_1(S,{z_0})$-action on $\tilde E$

$$(z,w) \rightarrow (\tilde \rho _{can}(\gamma) z,\tilde \rho (\gamma)  w)
\hskip 1cm,\hskip 1cm
\gamma \in \pi_1(S,{z_0}).
\eqno (1.3)$$
The quotient of $\tilde E$ by this action gives rise to a
vector bundle $E_{\tilde\rho}$ over $S$.
On $\tilde E$ we can consider the equation given
by $\tilde A=0$ (i.e. $\frac {dw}{   dz} =0$).
Its solutions are the constant functions.
Since this equation $\tilde A$ is invariant under
the action in (1.3), it descends to a linear ordinary
differential equation on $E_{\tilde\rho}$
which is holomorphic over  $S$. The construction gives
directly that the monodromy transformation of this
equation is the given representation
$\tilde \rho $.
The graphs of the local solutions
to $(1.1)$ form a holomorphic
foliation $\cF_{\tilde \rho}$ in $E_{\tilde \rho}$.

\subsection{ The Riccati Equation\label{ss.Riccati}}
\vskip 2mm
 Riccati equations may be obtained from
a linear ordinary differential equation
as $(1.1)$ or $(1.2)$  by projectivising the linear variables
of the vector bundle $E_{\tilde \rho}$ over the Riemann surface $S$.
Denoting $\zeta_j:=\frac{w_j}{  w_1}$ with $j=2,\ldots,n$,
the Riccati equation associated to (1.1) 
in affine coordinates takes the form
of a quadratic polynomial in
$\zeta_2,\ldots,\zeta_n$ with rational
coefficients in $z$:

{\small
$$\left(\begin{array}{c}
                                 \frac{d\zeta_2}{  dz} \\
                                  \cdots  \\
					\frac{d\zeta_n}{ dz}\\	
                    \end{array}\right)
=
\left(\begin{array}{c}
                                a_{21}\\
					\cdots\\
                                 a_{n1}\\
                    \end{array}\right)
+
\left(\begin{array}{ccc}
                                a_{22} - a_{11} & a_{23} & \cdots\\
                                  a_{32} & a_{33}-a_{11} & \cdots \\
					\cdots & \cdots &  a_{nn}-a_{11}\\

                    \end{array}\right) 
\left(\begin{array}{c}
                              \zeta_2 \\ \cdots \\ \zeta_n \\
                      \end{array}\right)
- (a_{12}\zeta_2 +\cdots + a_{1n}\zeta_n)
\left(\begin{array}{c}
\zeta_2 \\
\ldots \\
\zeta_n\\
\end{array}\right)     \eqno (1.4)
$$}where $A=(a_{ij}(z))$ is the matrix of
rational functions
in (1.1).
Similarly, given a representation $\tilde\rho$
as in $(1.2)$ we may also construct
from the projectivised representation $\rho$
in $(1)$
its suspension $M_\rho=Proj(E_{\tilde\rho})$
which gives a manifold which is a
$\CC P^{n-1}$ bundle over $S$ 
with a flat connection. The set of flat
sections form a foliation $\cF_\rho$ of $M_\rho$
which is the projectivisation of the foliation
$\cF_{\tilde\rho}$ in $E_{\tilde \rho}$.
The foliations so constructed,
will be called {\bf Riccati foliations}.

\section{The Foliated Geodesic Flow on
Linear and Riccati Equations\label{s.geodesic}}

\subsection{ The Geodesic Flow on Finite Hyperbolic Riemann Surfaces\label{ss.geosurface}}

 We say that $S$ is a finite hyperbolic Riemann surface if
$S$
is conformally equivalent
to  $\bar S - \{p_1,\ldots,p_r\}$, where $\bar S$ is a compact
Riemann surface of genus $g$ and
$g>1$ or $g=1$ with $r\geq1$ or $g=0$ with $r\geq3$. In such a case
$S$ has as a universal covering space
the Poincar\'e upper half plane $\HH^+$,
with its complete metric of curvature -1
given by
 $ds =\frac{ \vert dz\vert}{  y}$.
We introduce on $S$ the
hyperbolic metric induced by
the Poincar\'e metric via
the universal covering map. For the measure associated to the hyperbolic metric, the surface $S$ has finite area.

\vskip 2mm
Let $T^1S$ be the unit tangent bundle of $S$.
The Liouville measure $dLiouv$ on $T^1S$ is the measure
obtained from the  hyperbolic area element in $S$
and Haar measure $d\theta$ on unit vectors,
normalised so as to have volume 1.
The geodesic flow
$$\varphi:T^1S \times \RR \rightarrow T^1S
\hskip 1cm \varphi_t := \varphi(\ ,t)
\eqno (2.1)$$
is obtained by flowing along the geodesics (see [14] p. 209).
The geodesic flow leaves invariant the Liouville measure.

\vskip 2mm
\noindent
{\bf Theorem 2.1 (Hopf-Birkhoff).}  ([14] p. 217, 136)
{\it Let $S$ be a finite hyperbolic Riemann surface,
then the Liouville measure
is ergodic with respect to the geodesic flow
and the generic geodesic of $S$ is statistically
distributed in $T^1S$
according to the Liouville measure:
For all Liouville integrable functions $h$ on $T^1S$
and for almost any
$v_z \in T^1S$ with respect to the Liouville measure}

$$\lim_{t \rightarrow \infty}\frac {1}{  t} \int_0^th(\varphi(v_z,t))dt
= \int_{T^1S}hdLiouv$$

\subsection{The Foliated Geodesic Flows\label{geoRiccati}}

Let
$S$ be a finite  hyperbolic Riemann surface,
and $\tilde \rho$ and $\rho$ representation
as in $(1)$ and let $\cF_{\tilde\rho}$ and
$\cF_\rho$ be the foliations constructed in section 1.
If $\cL$ is a leaf of the foliation $\cF_{\tilde\rho}$
or $\cF_\rho$,
then the projection map $p:{\cL} \rightarrow S$
is a covering map, and hence the pull back of the
Poincar\'e metric of $S$ induces a metric to
the leaves of ${\cF}$, which coincides with
the  Poincar\'e metric of each leaf
$\cL$ of  ${\cF}$. This is
the Poincar\'e metric of the foliations
${\cF}_{\tilde \rho}$ or $\cF_\rho$.

\vskip 2mm
Let $T^1{\cF}_{\tilde \rho}$ 
be the manifolds
formed by those tangent vectors to $E_{\tilde\rho}$
and $M_\rho$ which are tangent to
 ${\cF}_{\tilde \rho} $ and $\cF_\rho$
 and are
of unit length with respect to the
 Poincar\'e metrics of the foliations.
The derivative of the projection map $E_{\tilde \rho},
M_\rho \rightarrow S$  induces
the commutative diagram

$$
\begin{array}{ccc}
T^1 {\cF}_{\tilde\rho} & {\overset
q \to } & E_{\tilde\rho} \\
 \downarrow&& \downarrow \\
T^1 {S} & {\overset q \to } & S \\
\end{array}
\hskip 2cm
\begin{array}{ccc}
T^1 {\cF}_\rho & \rightarrow & M_\rho \\
 \downarrow&& \downarrow \\
T^1 {S} & \rightarrow & S \\
\end{array} 
$$
The foliated geodesic
 flows $\tilde \Phi$ and $\Phi$
 are defined by following geodesics along the leaves
and is compatible with the geodesic flow
$\varphi$ on $S$, giving
 rise to the commutative diagram

$$
\begin{matrix}
\tilde\Phi:&
T^1{{\cF}_{\tilde\rho}}&
\times& 
\RR&
\rightarrow&
T^1{{\cF}_{\tilde\rho}}\cr
&
\downarrow&
&
\downarrow&
&
\downarrow
\cr
\varphi:&
T^1S&
 \times & \RR&
\rightarrow&
T^1S\cr
\end{matrix}
\hskip 10mm,
\hskip 10mm
\begin{matrix}
\Phi:&
T^1{{\cF}_{\rho}}&
\times & \RR&
\rightarrow&
T^1{{\cF}_{\rho}}\cr
&
\downarrow&
&
\downarrow&
&
\downarrow
\cr
\varphi:&
T^1S&
 \times& \RR&
\rightarrow&
T^1S\cr
\end{matrix}.
\eqno (2.2)
$$

\vskip 2mm
For any $v \in T^1S$ and $t \in \RR$, the flow $
\tilde\Phi_t:=\tilde\Phi(\ ,t)$
induces a
linear isomorphism
$$\tilde A(v,t):=\tilde \Phi(v,\ ,t)|_{E_{\tilde \rho,v}}:
E_{\tilde \rho,v} \rightarrow E_{\tilde\rho,\varphi(v,t)}
\eqno (2.3)$$
between the $\CC^{n}$-fibres.
After a measurable trivialisation
of the bundles by choosing measurably an othonormal basis of the fibers,
the foliated geodesic flows
may be seen as
measurable multiplicative cocycles
over the geodesic flow
on $T^1S$:
$$\tilde A: T^1S \times \RR \rightarrow GL(n,\CC)
\hskip 5mm,\hskip 5mm
\tilde A(v,t_1+t_2) = \tilde A(\varphi(v,t_1),t_2)  \tilde A(v,t_1)
\hskip 5mm,\hskip 5mm t_1,t_2\in\RR.
\eqno(2.4)
$$
Moreover the usual operator norm in $GL(n,\CC)$ coincides with the operator
norm of $(2.3)$ as Hermitian spaces with the metrics induced from the
fibre bundle metric.

\section{SRB-measures for Riccati Equations}
\subsection{SRB-measures\label{ss.srbdefi}}

Let $M$ be a differentiable manifold. The Lebesgue measure
class is the set of measures whose restriction on any chart $U$
has a smooth strictly-positive Radom-Nikodyn derivative with
respect to $dx_1\wedge dx_2\cdots\wedge dx_n$ where the $x_i$ are
coordinates on $U$.
 A set $E\subset M$ has {\em zero Lebesgue measure} if there is
 a measure $\mu$ in the Lebesgue class such that $\mu(E)=0$.

Let $X$ be a complete vector field on the manifold $M$, and denote
by $\varphi_t$ its flow.
A probability measure $\mu$ on $M$ is {\em invariant by $X$} if for
any $t\in\RR$ one has 
$\varphi_{t*}(\mu) = \mu.$
The {\bf basin}
 $B(\mu)$ of an $X-$invariant probability $\mu$ is the set of points
 $x \in M$ such that the
positive time average along its orbit of any continuous function $h \colon M\to\RR$ with compact support coincides with the integral of the function by $\mu$. In formula:
$$ \lim_{T\to+\infty}\frac1T\int_0^T h (\varphi_t(x))dt 
=\int_M h  d\mu$$.

\begin{defi} An $X-$invariant probability
measure in $M$ is an
{\bf SRB-measure} if its basin has non-zero Lebesgue measure
in $M$. 
\end{defi}

\subsection{Key Idea to Build SRB-measures for Riccati Equations}

\vskip 2mm
Let $S$ be a finite hyperbolic Riemann surface and
$\tilde \rho$ and $\rho$ representations as in $(1)$ and
$\cF_{\tilde \rho}$ and $\cF_\rho$ the foliations in
$E_{\tilde\rho}$ and $M_\rho$ 
constructed in section 2.
Consider a continuous Hermitian metric $|\cdot|_x$ 
on the fiber $E_{\tilde\rho,x}$ 
of $E_{\tilde\rho}$ and for each point $x\in S$ we endow
the corresponding Fubini-Study (Hermitian)
metric $ | \cdot | _x$ on $M_{\rho,x}= Proj(E_{\tilde\rho,x})$.
The bundles $ q^*E_{\tilde \rho} \simeq
T^1 \cF_{\tilde \rho}
 $ and $ q^*M_\rho  \simeq T^1 \cF_\rho $
over $T^1S$ are endowed in a natural way with the induced  Hermitian or  
Fubini-Study metric, respectively. 

   \begin{defi}
 Under the above setting,
 assume that the vector bundle $E \colon=
T^1\cF_{\tilde \rho} 
\to T^1S$ admits a measurable  splitting
$E_v = F_v\oplus G_v$ , defined for $v$ in a subset $\cA$ of $T^1S$,
and verifying the following hypothesis:
\begin{enumerate}
\item $\cA$ has total Lebesgue measure in $T^1S$;
\item $\cA$ is invariant by the geodesic flow $\varphi$;
\item the splitting is invariant by the foliated geodesic flow
$\tilde \Phi$:
for every $t\in\RR$ and every $v\in \cA$, $$F_{\varphi(v,t)}=
\tilde \Phi(F_v,t) \quad\mbox{and}\quad
G_{\varphi(v,t)}=\tilde\Phi(G_v,t);$$
\item $dim(F_v)=1;$
\item for any $v\in\cA$, for any compact set $K\subset T^1S$
and for any sequence $(t_n)_{n\in\NN}$ of times
such that $\varphi(v,t_n)\in K$ and
$\lim_{n\to\infty}t_n=+\infty$, one has:
$$\lim_{n\to \infty} \frac{ |\tilde \Phi(w_1,t_n) | _{\varphi(v,t_n)}}{ |
\tilde \Phi(w_2,t_n) | _{\varphi(v,t_n)}} = \infty,
\quad\mbox{for all non-zero}\quad w_1\in F_v,
\quad\mbox{and}\quad w_2\in G_v.
$$
\end{enumerate}
Under the above  hypothesis denote by
$\sigma^+\colon\cA\subset T^1S\to T^1\cF_\rho$ the
 mesurable section
defined by letting $\sigma^+(v)$ be the point of
$Proj(E_v)$ corresponding to the line $F_v$.
A  section $\sigma^+$ verifying the above hypothesis
is called a
{\bf  section of largest expansion}.
\end{defi}
Similarly, one defines
the {\bf  section
$\sigma^-$ of largest contraction}
by requiring

$$\lim_{n\to \infty} \frac{ |\tilde \Phi(w_1,t_n) | _{\varphi(v,t_n)}}{ |
\tilde \Phi(w_2,t_n) | _{\varphi(v,t_n)}} = \infty,
\quad\mbox{for all non-zero}\quad w_1\in F_v,
\quad\mbox{and}\quad w_2\in G_v.
\eqno (3.1)
$$
with $\lim_{n\to\infty}t_n=-\infty$
where we are imposing the condition that
the measurable sub-bundle
  $F$ is 1 dimensional (i.e. greatest expansion for negative times).

\vskip 2mm
\noindent
{\bf Proof of Theorem 1:}
$\sigma^+$ induces an isomorphism of the measure $dLiouv$ and $\mu^+
= \sigma^+_*dLiouv$,
so that the invariance and the ergodicity
of $\mu^+$ follow from those of $dLiouv$ and of $\sigma^+$.

\vskip 2mm
Let $h\colon T^1\cF_\rho\to \RR$ be a continuous function with compact support,
and denote by $K$ the projection of this compact set on $T^1S$. The function $h\circ \sigma^+\colon T^1S\to \RR$ is  measurable and bounded, so it belongs in $\cL^1(dLiouv)$.
As the Liouville measure is a 
$\varphi$ ergodic probability on $T^1S$,
there is an invariant set $Y_h\subset T^1S$ of total Lebesgue
measure such that, for $v\in Y_h$, the average

$$\frac1T\int_0^{T} h\circ\sigma^+
(\varphi(v,t))dt \ \rightarrow\ 
\int_{T^1S} h\circ\sigma^+ dLiouv=\int_{T^1\cF_\rho} h d\mu^+.
\eqno (3.2)$$

For each $v\in Y_h$ we denote by $\cY_h(v)$ the set of points in the fiber $y\in Proj(E_v)$ corresponding to a line of $E_v\setminus G_v$. We denote by $\cY_h$ the union
$\cY_h=\bigcup_{v\in Y_h} \cY_h(v) \subset M_\rho$.The set $\cY_h$ is invariant by
$\Phi$ because $Y_h$ is invariant by $\varphi$ and the bundle $G$ is
$\tilde\Phi-$invariant. By Fubini's theorem, the set
$\cY_h$ has total Lebesgue measure in $M_\rho$.
\begin{claim} For every $w\in \cY_h$, the average 
$\frac1T\int_0^{T} h(\tilde\Phi(w,t))dt$ converges to 
$\int_{T^1\cF_\rho} h d\mu^+$
\end{claim}
Before proving the claim let us show that this concludes the proof of Theorem 1:
There is a countable family $h_i, \quad i\in\NN$ of continuous functions with compact support which
is dense (for the uniform topology) in the set of all continuous 
functions of $T^1\cF$ with compact support. Look now at the set $\cY=\bigcap_0^\infty \cY_{h_j}$ : It is invariant by $\Phi$, has total Lebesgue measure, and is contained in the basin of $\mu^+$
by the claim. This proves  Theorem 1.

\vskip 2mm
Now we prove the claim:
Let $w\in \cY_h(v)$, for some $v\in Y_h$, and denote $w_0=\sigma^+(v)$.
As the section $\sigma^+$ is invariant by the foliated geodesic flow,
for any $t$, $\Phi(w_0,t)= \sigma^+(\varphi(v,t))$;
so  for any $T\in \RR$ the averages
$\frac1T\int_0^{T} h \circ \Phi(w_0,t)dt$ and $\frac1T\int_0^{T} h\circ\sigma^+(\varphi(v,t))dt$ are equal and we get 
by $(3.2)$
$$\lim_{T\to\infty}\frac1T\int_0^{T} h(\Phi(w_0,t))dt=
\int_{T^1\cF_\rho} hd\mu^+.$$

Consider a non-zero vector $\tilde w$ in the linear space $E_v$ in the line corresponding to
$w$. We can write in a unique way
$\tilde w= \tilde w_0+\tilde w_1$ where $\tilde w_0\in F_v$ and $\tilde w_1\in G_v$.
Notice that $\tilde w_0\neq 0$ projects on $w_0 \in Proj(E_v)$.
By hypothesis 5 in Definition 3.2,
when $t\in \RR$ is very large, either $\varphi_t(v)\notin K$
 or $\frac{ |\tilde \Phi( w_0,t_n) | _{\varphi(v,t_n)}}{|
 \tilde \Phi
 ( w_1,t_n) | _{\varphi(v,t_n)}}$ is very large and so the
 distance (for the Fubini-Study metrics) between
 $\Phi(w,t)$ and $\Phi(w_0,t)$ is very small, and goes to zero. Now we decompose the averages
$\frac1T\int_0^{T} h(\varphi_t(w))dt$ in two parts, one corresponding to the times $t$ such that
$\varphi(v,t)\notin K$, and the other to the times such that $\varphi(v,t)\in K$. 
The first part is uniformly zero (for both $w$ and $w_0$).
Moreover for large $t$ such that
$\varphi(v,t)\in K$, the difference $ h(\Phi(w_0,t))-h(\Phi(w,t))$
goes to zero.
So the averages of $h$ along the orbits of $w$ and $w_0$ converge to
the same limit, which is 
$\int_{T^1\cF_\rho} hd\mu^+$.
\qed

\begin{rema}\label{r.bigger}
\begin{enumerate}{\em
\item The existence of a  section of largest expansion does
not depend of the choice of the continuous Hermitian metrics on the fibers.
 
\item Theorem   does not use our specific hypotheses
(2-dimensional basis, geodesic flow, holomorphic foliation).
One has:

\vskip 2mm
\noindent
{\bf Theorem $1^\prime$:} 
{\it Let $B$ be a manifold and $\varphi$ a flow on $B$ admiting an
ergodic invariant probability $\lambda$ which is absolutely
continuous (with strictly positive density) with respect to
Lebesgue measure. Let $\tilde\rho\colon\pi_1(B) \to GL(n, \CC)$
be a representation, $(E_{\tilde\rho},\tilde\cF_{\tilde\rho})$ be the
vector bundle endowed with the suspension foliation, and
$M_\rho=(Proj(E_{\tilde\rho}), \cF_\rho)$
the suspension of the corresponding representation
$\rho\colon\pi_1(B)\to PGL(n, \CC)$.
Let $\Phi$ be the lift of 
the flow $\varphi$ to the leaves of $\cF_\rho$.  If
the bundle $E_{\tilde\rho}$ admits a section
$\sigma^+$ of largest expansion
then $\sigma^+_*(\lambda)$ is an SRB-measure of the flow $\Phi$,
whose basin has total Lebesgue measure in $M_\rho$.}

\vskip 2mm 
\item \label{r.sym} The geodesic flow (and the foliated geodesic flow) have a symmetry: 
denote by $I$ the involution map on the unit tangent bundle sending
each vector $v$ to $-v$ and $\tilde I$ the involution
$\tilde I(w_v) = -w_v$ on $T^1\cF_{\tilde\rho}$.
Then $I$ is a conjugation between the geodesic flow
and its inverse $I\circ\varphi_t\circ I=\varphi_{-t}$.
This shows that $\sigma^-=\tilde I\circ \sigma^+\circ I$
is a section
of largest expansion for the negative geodesic flow,
and $\mu^-=\sigma^-_*(dLiouv)$ will be an SRB-measure for the negative orbits of the geodesic flow.
Then Lebesgue almost every orbit in $T^1\cF$ has negative average
converging to $\mu^-$ and positive average converging to $\mu^+$.
}\end{enumerate}
\end{rema}

\begin{prop}
Let $E_{\tilde\rho}=F \oplus G$ be a $\tilde \Phi$-invariant
measurable splitting
giving rise to a section of largest expansion
$\sigma^+:=proj(F)$,
then the decomposition is measurably unique
(i.e. over a set of full Liouville measure
in $T^1S$).
\end{prop}

\begin{demo}
Let 
$E_{\tilde\rho}=F_1 \oplus G_1$ be a
$\tilde \Phi$-invariant  measurable splitting
giving rise to a section of largest expansion,
$\sigma_1^+:=proj(F_1)$. The line bundle  $F_1$ is not contained in $G$,
for if it were contained,
then the order of growth of $\sigma^+$ would be larger
than the order of growth of $\sigma_1^+$.
But then $G_1$ would not be a subset of $G$
and any initial condition in $G_1-G$
has the same order of growth than $\sigma^+$,
which is larger than the order of growth
of sections in $G$, like $\sigma_1^+$,
contradicting that the order of growth
of $\sigma_1^+$ is larger than the
order of growth of any section in $G_1$.

\vskip 2mm
Assume that $F \neq F_1$.
For $\varepsilon>0$  define the subset
$$H_\varepsilon:=\{v \in T^1S\ / \
dist(\sigma^+(v),G_v)>\varepsilon\ , \
dist(\sigma_1^+(v),G_v)>\varepsilon\ , \
dist(\sigma^+(v),\sigma_1^+(v))>\varepsilon\  \}
$$
where the distances are measured in the Fubini-Study metrics
of $Proj(E_v)$. For small $\varepsilon$ the set
$H_\varepsilon$ will have positive Liouville measure.
But since the Liouville measure is ergodic,
almost all points in $H_\varepsilon$ are recurrent.
But this cannot be, since both $\sigma^+$ and
$\sigma_1^+$ are invariant and as time increases
the component in $F_v$ grows much more than
the component on $G_v$ so that in $Proj(E_v)$
the sections $\sigma^+$ and
$\sigma_1^+$
are getting closer
which contradicts the condition
$dist(\sigma^+(v),\sigma_1^+(v))>\varepsilon$.
Hence we must have $F=F_1$ (Liouville almost everywhere),
as well as $\sigma^+=\sigma_1^+$. Now $G$ is uniquely determined
by $\sigma^+$, since any section outside $G$ has the same
order of growth as $\sigma^+$, and those on $G$ have
smaller order of growth.
\end{demo}

\section{Using Oseledec's Theorem\label{ss.oseledec}}

\subsection{ A Corollary of Oseledec's Theorem}

Let $$f:B\to B \hskip 2cm,\hskip 2cm
A\colon B\to GL(n, \CC)$$ be  measurable  maps.
For any $n\in\NN$ and any $x\in B$ we denote
$$A^n(x)= A(f^{n-1}(x)) \cdots A(f(x)) A(x) \hbox{ and }
A^{-n}(x)= [A^n(f^{-n}(x))]^{-1}.$$
One says that the family
$\{A^n\}$ form a multiplicative cocycle over $f$.

\begin{defi}
A point $x\in B$ has Lyapunov exponents
for the multiplicative cocycle $\{A^n\}$ over $f$ if
there exists $0<k\leq n$ and for all $i\in \{1,\dots,k\}$
there is $\lambda_i\in \RR$ and a subspace $F_i$ of $\RR^n$ such that:
\begin{enumerate}
\item $\RR^n = \bigoplus_i F_i$
\item For any $i$ and any non zero vector $v\in F_i$ one has 
$$
\lim_{n\to\pm\infty}\frac1n \log ( | A^n(v) | ) = \pm \lambda_i
$$
\end{enumerate}  
\end{defi}

\vskip 2mm
\noindent  
{\bf Oseledec's Multiplicative Ergodic Theorem
([14],p.666-667):}
Let $f: B\to B$ be an
invertible measurable transformation,
 $\mu$  an $f-$invariant probability measure
 and $A$ a mesurable multiplicative cocycle
over $f$.
Assume that the functions $\log^+\|A\|$ and $\log^+\|A^{-1}\|$
belong to  $\cL^1(\mu)$. Then the set of points for which the
Lyapunov exponents  of $A$ are well defined has $\mu$-measure $1$.
If $\mu$ is ergodic the Lyapunov exponents are independent of
the point in a set of 
total $\mu-$measure. 
  
 \vskip 2mm 

The Lyapunov exponents and the
Lyapunov spaces above depend measurably of $x\in B$ on a
set of $\mu-$total measure (see [14] p.666-667).
When the measure $\mu$ in Oseledec's Theorem is ergodic, we
can then speak of the {\em Lyapunov
exponents of the measure $\mu$}.

\vskip 2mm
We want  to use Oseledec's Theorem for
flows  when the base manifold is non-compact.
Let $\varphi$ be a complete flow on the manifold $B$,
 $\pi\colon E\to B$  a vector bundle over $B$ and
$\tilde \Phi$ be a flow on $E$ inducing a
multiplicative cocycle
as in $(2.4)$ over $\varphi$.

\begin{defi} We say that the Lyapunov exponents of $\tilde\Phi$ are well
defined at a point $v\in B$ if there is a continuous Euclidean or
Hermitian metric on the bundle $E$,
a finite sequence $\lambda_1<\cdots<\lambda_k$ and a
$\tilde\Phi$-invariant
splitting $E(v)=F_1(v)\oplus\cdots \oplus F_k(v)$
such that, for any non zero vector $w\in F_i(v)$,
any compact
$K\subset B$ and any sequence
$\{t_n\}_{n \in \ZZ }$
with $lim_{n \rightarrow \pm \infty} t_n = \pm \infty$ and
 $\varphi(v,t_n)\in K$ one has:
$$\lim_{ n  \rightarrow  \pm \infty} \frac1{t_{ n}}
\log( | \tilde\Phi(w, t_{ n}) | )=\pm\lambda_i.$$
\label{d.noncompact}
\end{defi}

The existence and the value of the Lyapunov exponents does
not depend of the continuous metric on the  vector bundle $E$;
moreover we can allow  the metric  to be discontinuous if  the
change of metric to a continuous reference metric is
bounded on compact sets of the basis $B$.

\begin{lemm} With the notation above the Lyapunov
exponents of  $v\in B$ for the 
 flow $\tilde\Phi$
are well defined if and only if they are well defined for the
multiplicative cocycle $\{\tilde A_1^n\}$ over $\varphi_1$
defined by the diffeomorphism $\tilde\Phi_1$.
Moreover the Lyapunov exponents and spaces  are
equal for the flow and the diffeomorphism.
\label{l.noncompact}
\end{lemm}
\begin{demo} One direction is clear, so we will assume that the
diffeomorphism $\Phi_1$ has Lyapunov exponents on $v$.
As the flow $\varphi$ is complete, for any compact set
$K\subset B$ the union 
$K_1=\bigcup_{t\in [-1,1]} \varphi(K,t)$ is compact.
Moreover for each $t_n$ such that $\varphi(v,t_n)\in K$,
let $T_n$ be the 
integer part of $t_n$, then $\varphi_1^{t_n-T_n}(v)\in K_1$.
We conclude the proof noticing that
$$\tilde A(v,t_n)= \tilde A(\varphi_1^{t_N-T_n}(v),T_n)
\tilde A(v,t_n-T_n)$$ and that the norm of
$\tilde A(*,T_n)$ is uniformly bounded
over $K_1$ independently of $t_n-T_n\in[0,1]$.
\end{demo}

\begin{defi} Let $\mu$ be a $\varphi-$invariant probability on $B$.
We say that the flow $\tilde\Phi$  defining a measurable
multiplicative cocycle $(2.4)$
is $\mu-$integrable if there is a
continuous norm $ | \cdot | $ on the vector bundle $E$
such that the functions $\log^+\|\tilde A_1\|$ and
$\log^+\|\tilde A_{-1}\|$ belong to $\cL^1(\mu)$,
where $\|\ \|$ is the operator norm on the normed vector spaces.
\end{defi}

The condition of integrability of the norm of the
multiplicative cocycle is always verified if the manifold $B$ is compact.

\vskip 2mm
  \noindent
 {\bf Proof of Corollary 2:}
 Consider $f=\varphi_1$, the time 1 
of the geodesic flow on $T^1S$,
and let $\tilde A(v)\colon E_v\to E_{f(v)}$ the linear multiplicative cocycle
 induced on the
vector bundle $T^1{\cF}_{\tilde\rho}$ by $\tilde\rho$
in Oseledec's Theorem.
By hypothesis,
this multiplicative cocycle is integrable so that the Lyapunov exponent
of the multiplicative cocycle $\tilde A$ are well defined for a
Liouville total measure set by Lemma 4.3.
The Lyapunov exponents
and spaces depend
measurably of $v\in T^1S$ which are invariant respectively by
$\varphi$ and $\tilde\Phi$. As the Liouville measure is ergodic, the Lyapunov exponents are constant on a set of total Liouville measure.
This ends the proof of item 1.

\vskip 2mm
The proof of item 2 is a direct consequence of the symmetry of the flow
$\Phi$: $\tilde I\circ \Phi_t\circ \tilde I= \Phi_{-t}$ (see
item 3 in  remark 3.3).
With the hypothesis of item 3 the section $\sigma^+$
is clearly a  section of largest expansion so that item 3
is a direct consequence of Theorem 1.
\qed

 \vskip 2mm
 A direct corollary of Theorem $1^\prime$ and
Oseledec's Theorem is the following

\vskip 2mm
\noindent
{\bf Corollary $2^\prime$:}{\it Let $f$ be a diffeomorphism of a manifold $B$,
admitting an invariant ergodic probability  $\lambda$ in the
class of Lebesgue and let $E$ be an $n-$dimensional vector bundle
over the basis $B$ and $M$  the corresponding projective bundle.
Assume that $\tilde \Psi$ is a diffeomorphism of $E$ leaving invariant the
linear fibration, inducing linear maps on the fibers and whose
projection on $B$ is the diffeomorphism $f$.
We denote by $\Psi$ the induced diffeomorphism on $M$.

Let $U_i$ be a covering of $B$ by trivializing charts of the bundle $E$:
then writing $\Psi$ in these charts we get a multiplicative cocycle
$\tilde A\colon B\to GL(n,\CC)$.
Assume that $log^+\|\tilde A\|$ and $log^+\|\tilde A^{-1}\|$ belong to $\cL^1(\lambda)$
and that the largest Lyapunov exponent of the measure $\lambda$ for the
multiplicative cocycle $\tilde A$ corresponds to a 1 dimensional space.
Denote by $\sigma^+$ the corresponding measurable section defined
on a Lebegue total measure set of $B$ to $M$.

Then $\sigma_*^+(\lambda)$ is an SRB-measure for $\Psi$ and
its basin has total Lebesgue measure in $M$.}
 \qed

\subsection{ Proof of Theorem 3}

\begin{demo}  Due to Corollary 2 and the Remark 3.2,
the only thing that remains to be proved
is that, under the integrablity condition $(3)$,
if there is a section of largest expansion
then the largest Lyapunov exponent is positive and simple.

\vskip 2mm
We begin first with the case that $S$ is compact.
So assume that there is a section $\sigma^+$ of largest expansion
providing a measurable decomposition $E_{\tilde \rho}=F\oplus G$,
$\sigma^+:=Proj(F)$
and let $\lambda_i$ and $F_i$ be
the Lyapunov exponents and spaces as in Corollary 2.
We have $F \subset F_k$,
corresponding to the greatest eigenvalue $\lambda_k$,
and denote by $H$ the measurable
bundle $F_k \cap G$ of dimension $n_k-1$.
Assume that the dimension $n_k$ of $F_k$ is at least 2,
and we will argue to obtain a contradiction to this assumption.

\vskip 2mm
Since the foliated geodesic flow leaves invariant the measurable
bundle $F_k$, after a measurable trivialisation we will obtain
a measurable cocycle

$$B:T^1S \times \RR
 \rightarrow GL(n_k,\CC)$$
which carries the information of how initial
conditions are transformed into final conditions,
when starting from the point $v \in T^1S\ , \ w \in F_{k,v}$, and
flowing a time $t$ along the geodesic.

\vskip 2mm
Recall that we have introduced a Hermitian metric on the
bundle $E_{\tilde{\rho}}$, by pull back in the bundle
$q^*E_{\tilde \rho} = T^1\cF_{\tilde \rho}$ and
by restriction into the bundle $F_k$.
Recall also that if we have a $\CC$-linear map $L$ between
Hermitian spaces, the determinant
$det(L,W)$ of $L$ on a subspace $W$
is by definition the quotient of the volumes of the
paralelograms determined by $Lw_1,\ldots,Lw_m,iLw_1,\ldots,iLw_m$
and $w_1,\ldots w_m,iw_1,\ldots,iw_m$ corresponding
to any $\CC$-basis $w_1,\ldots,w_m$ of $W$.
Define
$$\Delta^m:T^1S \rightarrow \RR \hskip 1cm,\hskip 1cm
\Delta^m(v):= \frac
{det(B(v,m),F_v)^{n_{k}-1}}{det(B(v,m),H_{v})}$$
and note that the cocycle condition $(2.4)$ for $B$
and the $\tilde \Phi$-invariance of $H$ and $F$ gives the
multiplicative condition
$$\Delta^m(v) = \Delta(\varphi(v,m-1)) \Delta(\varphi(v,m-2))
\cdots \Delta(v)
\hskip 5mm , \hskip 5mm \Delta:=\Delta^1.
\eqno(4.1)$$
The volume in $H$ has
exponential rate of growth  $(n_k-1)\lambda_k$,
since it is the Lyapunov exponent of $\Lambda^{n_k-1}H$.
The exponential rate of growth of $F$ is $\lambda_k$, hence

$$\int_{T^1S}\log(\Delta)dLiouv = (n_k-1)\lambda_k -
  \lambda_k - \ldots - \lambda_k  =0.
  \eqno (4.2)$$

\vskip 2mm
Now we need the following
corollary of a general statement from Ergodic Theory,
(see [15], Corollary 1.6.10):

\begin{coro} Let $\varphi:B \rightarrow B$ be a measurable transformation
preserving a probability measure $\nu$ in $B$, and $g:B \rightarrow \RR$
a $\nu$-integrable function such that
$\lim_{n \rightarrow \infty} \sum_{j=0}^n (g \circ \varphi^j) = \infty$
at $\nu$-almost every point, then
$\int_B g d\nu >0$.
\end{coro}

\begin{demo} 
 Consider the set
$$A:=\{v \in T^1S \ / \  \sum_{j=0}^\ell (g \circ \varphi^j)(v)>0,
\ \forall \ell\geq 0 \},$$
and for $v \in A$ let
$$S_*g(v):= \inf_\ell \{
  \sum_{j=0}^\ell (g \circ \varphi^j)(v)
\}.$$
$A$ has a strictly positive $\nu$ measure since almost any orbit will
have a point in $A$,
and
$S_*g$ is a measurable function on $A$ which is strictly positive. By
 Corollary 1.6.10 in [15] we have

$$\int_B g d\nu = \int_A S_*g d\nu,$$
but this last number is strictly positive,
since we are integrating a strictly positive function
over a set of positive measure.
\end{demo}

We want to apply the above Lemma to $(X,\nu) = (T^1S,dLiouv)$
and $g=\log\Delta$.
Note that the multiplicative relation $(4.1)$ implies
$$\sum_{j=0}^{m-1}\log\Delta(\varphi_j(v)) = \log \Delta^m(v)
\eqno (4.3)$$
The hypothesis on the growth  of the section $\sigma^+$
implies that
$\lim_{n\rightarrow \infty} \log\Delta^m(v) \rightarrow \infty$.
But using $(4.3)$
this is the hypothesis in the Lemma, so as a conclusion of
it
we obtain that
$$\int_{T^1S}\log(\Delta)dLiouv >0,$$
which contradicts $(4.2)$.
Hence $F_k$ has dimension 1, so that
the largest Lyapunov exponent is simple.

\vskip 2mm
Assume now that $S$ is not compact.
According to Lemma 4.3,
it is sufficient to consider
the integrability condition for the time 1 flow $\varphi_1$.
Let $K$ be a compact set
of positive Liouville measure
in
$T^1S$ and partition
$$K_m:=\{v \in K \ / \ \varphi^j(v) \notin K, j=1,\cdots,m-1,\ \varphi^m(j) \in K \}$$
according to the time of the first return to $K$.
Define the multiplicative cocycle
generated by $$
C:K \rightarrow GL(n,\CC)
\hskip 1cm,\hskip 1cm
C(v) :=  \tilde A_1^m(v) \hskip 1cm,\hskip 1cm v \in K_m$$
corresponding to the first return map to $K$.
Since
$$C(v) =
\tilde A_1(\varphi^{m-1}(v))
\ldots
\tilde A_1(\varphi(v))
\tilde A_1(v),               $$
we have
$$\log^+( \|C(v)\|) \leq
\log^+(\|\tilde A_1(\varphi^{m-1}(v))\|)+
\ldots
+\log^+(\| \tilde A_1(\varphi(v))\|)    +
\log^+(\|\tilde A_1(v)\|)),               $$
and hence on $K$ we obtain
$$\sum_{m=1}^\infty \int_{K_m}\log^+( \|C(v)\|) \leq
\sum_{m=1}^\infty[\log^+(\|\tilde A_1(\varphi^{m-1}(v))\|)+
\ldots
+\log^+(\| \tilde A_1(\varphi(v))\|)    +
\log^+(\tilde A_1(v))]\leq                $$
$$\leq  \int_{T^1S}\log^+( \|\tilde A_1(v)\|)$$
since the sets $$\varphi_j(K_m)
\hskip 1cm , \hskip 1cm j=0,\ldots, m-1, \hskip 1cm, \hskip 1cm
 m=1,\ldots$$
are disjoint. Hence the cocycle generated by $C$ is integrable,
and we may repeat the argument presented for the case
that $T^1S$ is compact.

\end{demo}

\section{Using Oseledec's Theorem in the Non-compact case}

The objective of this paragraph is  to prove
Theorem 4.
The proof of the parts "if" and "only if"
 are given by some estimates over the punctured disc $\DD^*$.
 As both  proofs are long,  we will treat them separately.
The common argument is the following estimate about the geodesic
flow of $\DD^*$.

\subsection{Estimates on the Geodesic Flow on a Punctured Disc}

Denote by $\DD^*$ the punctured disc
endowed with the usual complete metric of curvature $-1$,
that is, its universal cover is
the Poincar\'e half plane $\HH^+$ with covering group
generated by the translation $T(z)= z+1$ and define
$$D^*:= \frac{\{z\in \HH^+ \ / \ Im(z)>1\}  }{ (T^n)} \subset \DD^*
 \hskip 1cm, \hskip 1cm
S^1:=\partial D^*=
 \frac{\{z\in \HH^+ \ / \ Im(z)=1\} }{ (T^n)}
  \subset \DD^*.$$
$$\bar {D^*}:=
\frac{\{z\in \HH^+ \ / \ Im(z)\geq 1\} }{ (T^n)} \subset \DD^*
$$

A unit vector $u\in T^1D^*$ at a point $z\in D^*$ is called a
{\em radial vector} if $u \in \RR w
 \frac{\partial}{\partial w}$.
Note that for any non-radial vector $u\in T^1D^*$ the geodesic
$\gamma_u$ through
  $u$
  in $\bar{ D^*}$
  is a compact segment $\gamma_u$ whose extremities
  are on the circle $S^1$.
  We will denote the tangent vector of the geodesic $\gamma_u$
  on $S^1$ by $\alpha(u)$ (the incoming)
  and $\omega(u)$ (the outgoing), and
  let $t(u)$ be the lenght of $\gamma_u$.
The set of radial vectors has zero Lebesgue measure.
We will denote by $M$ the set of nonradial unit vectors on
$T^1\DD^*|_{\bar{ D^*}}$
and by $N$ the subset of $M$ over the circle
$S^1$. We denote $N^+$ the set of vectors
in $N$ pointing inside $D^*$
and by
$$\cA=\{(u,t), u\in N^+, t\in[0,t(u)]\}\subset N^+\times
[0,+\infty[.$$
The geodesic flow $\varphi$ on $T^1\DD^*$
induces a natural map $F\colon \cA\to M$ defined by $F(u,t)=\varphi(u,t)$.
The unit tangent bundle over $S^1$ admits natural coordinates :
If $u$ is a unit vector at $w$ we will denote $\theta(u)$
the argument of $w$, and $\eta(u)$ the angle between $u$
and the radial vector $-z\partial/\partial z$.
We denote by $\mu$ the measure on $\cA$ defined by $d\mu=
\cos(\eta)\cdot d\theta\wedge d\eta\wedge dt$

\begin{lemm} The Liouville measure on $T^1D^*$ is $F_*(d\mu)$
(up to a multiplicative constant).
\end{lemm}
\begin{demo} The measure $F_*^{-1}(dLiouv):=hd\theta\wedge d\eta\wedge dt$
for a certain function $h$. Since the Liouville measure is invariant
under the geodesic flow, and in $M$ the geodesic flow
has the expression $\frac{\partial
}{ \partial t}$, then $h$ is independent of $t$. Since the
Liouville measure is invariant under rotations in $\theta$
then $h$ is also independent of $\theta$. Hence $h$
is only a function of $\eta$.
To compute the value of $h$ it is enough to
compute for an arbitrary $\eta$ at a point in $N^+$.
We have $F_*(d\theta \wedge dt) = h(\eta) dArea$.
The variable $\theta$ is parametrized according to geodesic length
and since the angle between the vertical and the geodesic
at $Im(z)=1$ is $\eta$, we project the tangent vector to
the geodesic to the vertical direction to obtain the
weight $cos(\eta)$.
\end{demo}

We will denote by $\mu_0$ the measure on
$N^+$ defined by $d\mu_0=  d\theta\wedge d\eta$.

\begin{prop} Let  $\tilde A_t:T^1D^* \times \RR
\rightarrow GL(n,\CC)$ be a linear multiplicative cocycle over
the geodesic flow of $D^*$.
For every unit vector $u\in N^+$, we denote
$$B:N^+\rightarrow GL(n,\CC)
 \hskip 1cm,\hskip 1cm
 B(u) = \tilde A_{t(u)}(u)$$
 the matrix corresponding to  the geodesic $\gamma_u$
of length $t(u)$ going from
$\alpha(u)$ to $\beta(u)$. Then the two following sentences are equivalent:
\begin{enumerate}
\item There is a Hermitian metric $|\cdot|$ on the vector
bundle over $T^1D^*$
such that the multiplicative cocycle $\tilde A_1$ is integrable
for Liouville, that is
$$\int_{T^1D^*} \log^+\| \tilde A_{\pm 1}\| dLiouv<+\infty.
\eqno (5.1)
$$
\item The function $\log^+( \|B\|)$ belongs to $\cL^1(\mu_0)$, that is
$$\int_{N^+} \log^+(\| B(u)\| ) d\mu_0 <+\infty.
\eqno (5.2)$$
\end{enumerate}
\label{p.mu0}
\end{prop}
\begin{rema} {
$(5.2)$ does not depend of the choice of the continuous Euclidean metric :
Two continuous Hermitian metrics $|\cdot|_1$ and $|\cdot|_2$ on
the bundle over $T^1\DD^*|_{\partial D^*}$ are equivalent
because $\partial D^*$ is compact, so that the difference
$|\log(\| B(u)\| _1)| - |\log(\| B(u)\| _2)|$
is uniformly bounded on $N^+$.}
\end{rema}
\begin{demo} 
For every $u\in N^+$ set $t_u:=t(u)$, and
 divide the interval $[0, t_u]$ in
$$[0,1]\cup [1,2]\cup \cdots\cup [ E(t_u)-1,E(t_u)]\cup[E(t_u),t_u],$$
so that  if $u$ is a vector at a point $x \in\partial D^*$ one gets
on setting $\varphi:=\varphi_1$ the geodesic flow at time 1:

$$ B(u)= \tilde A_{t_u-E(t_u)}(\varphi^{E(t_u)}(u)\circ
\prod_0^{E(t_u)-1} \tilde A_1(\varphi^i(u))$$

So for any Hermitian norm $|\cdot |$  we get
$$ \|  B(u)\|  \leq \|  \tilde A_{t_u-E(t_u)}
(\varphi^{E(t_u)}(u))\|  \prod_0^{E(t_u)-1} \| \tilde A_1(\varphi^i(u))\| $$
So
$$ \log^+(\|  B(u)\| ) \leq \log^+ \|  \tilde A_{t_u-E(t_u)}(\varphi^{E(t_u)}(u))\|  +
\sum_0^{E(t_u)-1} \log^+\| \tilde A_1(\varphi^i(u))\| $$

Remark that $\log^+ \|  \tilde A_{t_u-E(t_u)}(\varphi^{E(t_u)}(u))\| $ is uniformly
bounded by a constant $K$ depending on $\tilde A$ and $|\cdot|$,
because  $t_u-E(t_u)\in[0,1[$ and $\varphi^{E(t_u)}(u) =
\varphi_{E(t_u)-t_u}(\varphi_{t_u}(u))$ remains in a compact set
(recall that $\varphi_{t_u}(u)\in\partial D^*$).
So we get that there is a constant $K_1$ such that for every
$u\in N^+$ one has

$$
\log^+(\|  B(u)\| ) \leq K_1 + \int_0^{t_u} \log^+\| \tilde A_1(\varphi_t(x))\|  dt
$$

Notice now that, for any $\varepsilon \in [0,1[$ there is $\delta>0$
such that if $\cos(\eta)\leq \varepsilon$ then $t_u\leq\delta$.
So it is equivalent that the function $\log^+(\| B\| )$ is integrable
for the measure $d\mu_0$ or for $\cos(\eta) d\theta\wedge d\eta$.

\vskip 2mm
Hence we obtain that if $\int_{N^+} \log^+(\| B\| ) d\mu_0 =+\infty$
then for any Riemannian metric $|\cdot|_2$ the function
$\log^+(\| \tilde A_1\| _2)$ is not Liouville integrable.
We have proven that
$ \mbox{item 1}\quad \Longrightarrow \quad \mbox{item 2}.$

\vskip 2mm
For the other implication, choose a continuous Riemannian
metric on the bundle  over $N$, assume
the integrability condition $(5.2)$ and let $v \in T^1\DD|_{\bar D^*}$.
If $v$ is a radial vector, then
push forward the metric over $\alpha(v)$ along the geodesic
using the flat structure of the bundle.
If $v$ is not a radial vector then push forward the
metric on $\alpha(v)$ on the first third of $\gamma_v$,
on the last third of the geodesic push forward
the metric on $\omega(u)$ and on the middle
third of $\gamma_u$ put the corresponding
convex combination of the metrics on $\alpha(u)$ and $\omega(u)$.
This produces a continuous metric on the bundle over
$T^1\DD|_{\bar D^*}$ such that $\|\tilde A_{\pm 1} \|$
does not expand except in the middle part,
and there it expands in a constant way.
Hence for this metric
the integral $(5.2)$ coincides with $(5.1)$.
\end{demo}

To use Proposition 5.2 we will need to
 estimate $\| B(u)\| , \quad u\in N^+$.
For that we will use the following estimate of $t_u$ and the estimate of the variation of the argument along the geodesic $\gamma_u$:

\begin{prop} 
\begin{enumerate}
\item There is a constant $T$ such that $t_u\in [-2\log |\eta|-T,-2\log|\eta| +T]$.
\item Denote by $a_u$ the variation of the argument along $\gamma_u$. Then 
$a_u= 2\frac{\cos \eta}{\sin\eta}$
\end{enumerate}
\label{p.tu}
\end{prop}
\begin{demo}
The easiest way is to look at the universal cover
$\HH$.
Recall that in this model the geodesic for the hyperbolic metric are circles or straight lines (for the Euclidean metric) orthogonal to the real line.
Let $u\in E^+_1$ at a point $ x\in \partial D^*$.
Denote by $u$ the corresponding vector at a point
$\tilde x\in \HH$, $Im(x)=1$, where $\tilde x$ is a lift of $x$.
The angle $\eta(u)$ is the angle between the vector and the vertical line.
Consider the geodesic $\tilde \gamma_u$ throught $u$.
The Euclidean radius $R_u$ of this circle verifies $1=|\sin(\eta)|\cdot R_u$.
Now denote by $\tilde y\neq \tilde x$ the intersection point
of $\tilde \gamma_u$ with the boundary $Im(z)=1$ of ${ {D}}^*$.
Then $a_u = \tilde y -\tilde x =
2\frac{\cos (\eta)}{\sin(\eta)}$.
So the second item of Proposition 5.4 is proved.

\vskip 2mm
To give an estimate of $t_u$ let us consider the following curve $\sigma_u$ joining the points
$\tilde x$ and $\tilde y$: $\tilde \sigma_u$ is the union of the vertical segment $\sigma^1_u$ joining $\tilde x=(\cR e(\tilde x),1)$ to $(\cR e(\tilde x), R_u)$ the horizontal segment $\sigma_2^u$ joining 
$(\cR e (\tilde x),R_u)$ to $(\cR e(\tilde y), R_u)$ and the vertical segment $\sigma^3_u$ joining $(\cR e(\tilde y),R_u)$ to $(\cR e(\tilde y), 1)=\tilde y$. 

\vskip 2mm
The hyperbolic length of the vertical segments is $\log (R_u)$.
The hyperbolic length of the horizontal segment is $\frac{|a_u|}{R_u}= 2cos(\eta)$. 
So we get:

$$\ell(\tilde \gamma_u) < \ell(\sigma_u)= -2\log(|\sin(\eta)|)+2\cos(\eta)$$

On the other hand, consider the point $z_u\in\gamma_u$ whose imaginary part is $R_u$. This point is the middle of the horizontal segment of $\sigma_u$. Denote by $\gamma_u^0$ the segment of $\gamma_u$ joining $\tilde x$ to $z_u$ and $\sigma_u^0$ the segment of $\sigma^2_u$ joining $zu$ to the point $(\cR e(\tilde y), R_u)$. The union of these 2 segments is a segment joining the two extremities of $\sigma_u^1$ which is a geodesic. So we get 

$$-\log(|\sin(\eta)|)<\ell(\gamma_u^0)+\ell(\sigma_u^0) =\frac12\ell(\tilde\gamma_u) +\cos(\eta).$$

So we get
$$t_u= \ell(\tilde \gamma_u)\in [-2\log(|\sin(\eta)|)-2\cos(\eta),-2\log(|\sin(\eta)|)+2\cos(\eta)] $$
So
$$t_u \in [-2\log(|\sin(\eta)|)-2,-2\log(|\sin(\eta)|)+2]$$

To conclude the first item it is enought to note that $|\log (|\eta|)-\log(|\sin(\eta)|)|$ is bounded for $\eta\in [-\pi/2,\pi/2]$.
\end{demo}

\subsection{The Parabolic Case}

\begin{prop}If for each $i$ all the eigenvalues of $\rho(\gamma_i)$ have  modulus $1$, then the multiplicative cocycle flow is integrable.
\label{p.if}
\end{prop}

As the function $\log^+ |\tilde A_1|$ is continuous, it is integrable for the Liouville  measure over every compact set of $T^1S$. So the problem is purely local, in the neighbourhood of the punctures of $S$.
So it is enough to look at a multiplicative cocycle $\tilde A_t$ over the
geodesic flow of the punctured disc $D^*$.
The proposition is a direct corollary of the following proposition:

\begin{prop} Let $B\in GL(n,\CC)$ be a matrix and $\cF_B$ be the corresponding suspension foliation over $\DD^*$ (as $B$ is isotopic to identity the foliation $\cF_B$ is on $\DD^*\times \CC^n$), and denote by $\tilde A_t$ the linear multiplicative cocycle  over the geodesic flow $\varphi$ of $\DD^*$ induced by $\cF_B$.
Assume that all the eigenvalues of $B$ have modulus equal to $1$.
Then the functions $\log^+(\|\tilde A_{\pm 1}\|)$ are in $\cL^1(dLiouv|_{D^*})$.
\label{p.parabolic}
\end{prop}
We begin the proof of Proposition 5.6 by the following remarks allowing us to reduce the proof to an easier case:

\begin{rema}
{\em
\begin{enumerate}
\item If two matrices $B_1$ and $B_2$ are conjugate then the corresponding cocycles are both integrable or both non-integrable.
\item If $B$ is a matrix on $\CC^k\times \CC^m$ leaving invariant $\CC^k\times\{0\}$ and $\{0\}\times\CC^m$, then the multiplicative cocycle induce by $B$ is integrable if and only if the cocycles induced by the restrictions of $B$ to $\CC^k\times \{0\}$ and $\{0\}\times\CC^m$ are both integrable.
\item As a consequence of  item 2, we can assume that $B$ is a
matrix which doesn't leave invariant any splitting of $\CC^n$ in a
direct sum of non-trivial subspaces.
In particular $B$ has a unique eigenvalue $\lambda_B$ and by hypothesis $|\lambda_B|=1$. Moreover two such matrices are conjugate: their Jordan form is
$$\left(\begin{array}{cccccc}
                                 \lambda_B&1&&\cdots&0&0\\
                                  0&\lambda_B&1&\cdots&&0\\
					\cdots&&&&&\cdots\\
					0&0&\cdots&0&\lambda_B&1\\
					0&0&\cdots&0&0&\lambda_B\\	
                    \end{array}\right)
$$ 

\end{enumerate}
}
\end{rema}

Using the remarks above, it is enough to prove Proposition 5.6 for the matrices $B_\theta$ define as follows.
Let 
$$A_\theta = \left(\begin{array}{cccccc}
                                 i\theta&1&0&\cdots&0&0 \\
                                  0&i\theta&1&0&\cdots&0  \\
					\cdots&&&&&\cdots\\
					0&0&\cdots&0&i\theta&1\\
					0&0&\cdots&0&0&i\theta\\	
                    \end{array}\right).
$$ We define $B_\theta=exp(A_\theta)$.  
Notice that 
$$ exp(t\cdot A_\theta)=e^{i t\theta}\left(\begin{array}{cccccc}
                                 1&t&&\cdots&\frac{t^{n-2}}{(n-2)!}&\frac{t^{n-1}}{(n-1)!} \\
                                  0&1&t&\cdots&&  \frac{t^{n-2}}{(n-2)!}\\
					\cdots&&&&&\cdots\\
					0&0&\cdots&0&1&t\\
					0&0&\cdots&0&0&1\\	
                    \end{array}\right)
$$

Consider the  holomorphic foliation 
defined by the linear equation 

$$\begin{pmatrix} \dot z \cr \dot w \end{pmatrix}  =
 \begin{pmatrix} i & 0 \cr 0 & A_\theta \end{pmatrix}  
\begin{pmatrix}  z \cr  w \end{pmatrix}  $$
on
$\DD^*\times \CC^n$ such that the holonomy map from $\{e^{-2\pi}\}\times \CC^n\to\{z\}\times\CC^n$
with $z \in S^1$ is $exp(arg(z)A_\theta)$. 
The monodromy of this foliation is $B_\theta = e^{2i\pi \theta} exp(2\pi A_0)$.

\begin{lemm} The multiplicative cocycle $\tilde 
A_t$ obtained by lifting the geodesic flow of $\DD^*$ on the leaves of $\tilde\cF_\theta$ is integrable over $T^1\DD|_{D^*}$.
\end{lemm}
\begin{demo} For any $u\in N^+$ one has $B(u)=A_{t_u}(u)=
exp(\frac{a_u}{2\pi} \cdot A_{\theta})$, so
that there is a constant $K$ such that
$\| B(u)\|< K(1+ a_u^{n-1})$, so that $\log^+\|B(u)\|$
is integrable if and only if $\log^+(|a_u|)$ is integrable for $\mu_0$.

By Proposition~\ref{p.tu} one has  $a_u= 2\cos(\eta)/\sin(\eta) $ so that $a_u<2/\eta$.
As $\int_{-1}^1 |\log(| 1/x|)|dx<+\infty$, we get easily that 
$\int_{N^+_1} \log^+(|a_u|) d\mu_0 <+\infty$, concluding the proof.
\end{demo}
\subsection{The Hyperbolic Case}

\begin{prop} If there is i such that the matrix $B=\rho(\gamma_i)$ has an eigenvalue with modulus different from $1$, then the multiplicative cocycle is not integrable.
\label{p.hyperbolic}
\end{prop}
If $B\in GL(n,\CC)$ has an eigenvalue with modulus different from  $1$,
we may suppose that its modulus is greater than $1$, since 
 the suspension of $B$ and $B^{-1}$ are isomorphic.
As in the parabolic case the proof of Proposition~\ref{p.hyperbolic} follows directly from a local argument in a neighbourhood of the puncture corresponding to $\gamma_i$.

\begin{prop} Let $B\in GL(n,\CC)$ having an eigenvalue $\lambda>1$ and $\cF_B$ the suspension folition on $D^*$. Then the multiplicative cocycle $\tilde A_t$ induced by $\cF_B$ over the geodesic flow $\varphi$ of $D^*$ is not integrable.
\label{p.localhyp}
\end{prop}
\begin{demo}  
We begin by an estimate of the norm of the multiplicative cocycle
 corresponding to the "in-out"  map :
\begin{lemm}There is a constant $K>0$ such that for any $u\in N^+$ one has: 
$$| \tilde A_{t_u}(u)| \geq K\cdot \lambda^{a_u/2}.$$
\end{lemm} 
So $\log^+|\tilde A_{t_u}(u)|\geq \log K + \frac{|a_u|}2 \log\lambda$. One deduces that $\log^+|\tilde A_{t_u}(u)|$ cannot be $\mu_0$-integrable if $|a_u|$ is not integrable.
By Proposition~\ref{p.tu} one knows that 
$a_u = 2\frac{\cos(\eta)}{\sin(\eta)}$ 
and this function is not integrable for $d\mu_0= d\eta\wedge d\theta$.
From Proposition~\ref{p.mu0} we get that the multiplicative cocycle  $\tilde A_1$ is not integrable for Liouville, finishing the proof the Proposition~\ref{p.localhyp}.
\end{demo}

\vskip 2mm
\noindent
{\bf Remark:} 
If $\rho:\pi_1(S) \rightarrow PGL(n,\CC)$ 
is a representation that does not admit a lifting to a representation in $GL(n,\CC)$ we may still define a flat bundle over $S$ 
but with fibres $\CC^n/\ZZ_n$ and transition coordinates in
$SL(n,\CC)/\ZZ_n\cdot Id$, and hence a foliation $\cF_{\tilde \rho}$
on this singular bundle, where $\ZZ_n$ is the group of $n$ roots of unity.
We may introduce a continuous Hermitian norm on this bundle 
(locally induced from a Hermitian norm in $\CC^n$ 
as well as choosing a trivialisation of the generator
of the discrete dynamics $\tilde A_1$,
and the statements and arguments given in the text
extend to this situation.

 \vskip 5mm
\section{Ping-pong and Schottky Monodromy  Representations}
\label{s.pingpong}

The ping-pong is a classical technique used to verify that a finitely
generated group of transformation of some space is a free group.
When the space is a metric space additional geometric information
on the ping-pong allows one to describe almost completely the
topological dynamics of this group of transformations.
 We will use this technique to describe the foliated geodesic
 flow associated to an
 injective  representation $\rho$ from $\pi_1(S)$ to a Schottky group $\Ga\subset PSL(2,\CC)$.

\subsection{The Ping-pong}

Let us first recall some basic properties and definitions on the ping-pong.

\begin{defi} Let $\cE$ be a set, $k>1$  and  for every $i\in\{1,\dots,k\}$ let $f_i\colon\cE\to\cE$ be a bijection. We say that {\em the group $\Ga\subset \mbox{Bij}(\cE)$ generated by $f_1,\dots,f_k$ is a ping-pong (for this system of generators)} if for every $i\in\{1,\dots,k\}$ there exist subsets $A_i$, $B_i$ of $\cE$ such that the following properties are verified:
\begin{itemize}  
\item The family $\{A_i, B_i, i\in\{1,\dots,k\}\}$ is a family of mutualy disjoint subsets of $\cE$,
\item for every $i\in\{1,\dots,k\}$ one has $f_i(\cE\setminus A_i)\subset B_i$.
\end{itemize}
\label{d.pingpong}
\end{defi}

Denote by $\FF_k$ the free group with $k$ generators $\{e_1,\dots,e_k\}$. The first result on the ping-pong is:

\begin{prop} If a group $\Ga\subset\mbox{{\em Bij}}(\cE)$ is a ping-pong group for the generators $f_1,\dots,f_k$ then the morphism $\varphi\colon\FF_k\to \Ga$ defined by $\varphi(e_i)=f_i,\quad i\in\{1,\dots,k\}$ is an isomorphism.
\end{prop}
\begin{demo} Let $i_1,\ldots,i_m \in \{1,\ldots,k\},
$ and $\varepsilon_j\in\{-1,1\}$ be
such that the word
$e_{i_1}^{\varepsilon_1}\cdots e_{i_m}^{\varepsilon_m}$ is
a reduced word in $\FF_k$. We have to prove that the bijection
$f= f_{i_m}^{\varepsilon_m}\circ\cdots\circ f_{i_i}^{\varepsilon_1}=
\varphi(e_{i_1}^{\varepsilon_1}\cdots e_{i_m}^{\varepsilon_m})$ is
different from identity. For instance assume
that $\varepsilon_1=1$. Then, using that the word is a reduced word,
one easily shows (by induction on $m$) that $f(\cE\setminus A_{i_1})$ is included in one of the sets $A_{i_m}$ or $B_{i_m}$. As $k>1$,  $f(\cE\setminus A_{i_1})$ is not included in one element of $\{A_i, B_i, i\in\{1,\dots,k\}\}$, so $f$ is not the identity.   \end{demo}

Assume now that $(\cE,d)$ is a  compact metric space,
the $f_i$ are homeomorphisms of $\cE$, every $A_i$, $B_i$ is compact,
and for each $i\in\{1,\dots,n\}$ the restrictions of $f_i$ and
$f^{-1}_i$ to $\cE\setminus A_i$ and $\cE\setminus B_i$, respectively,  are
contractions for the distance $d$: we will say that $(\cE,d,\{f_i\})$ is a {\em compact contracting ping-pong. }

\vskip 2mm
For any $g\in\{f_i, f_i^{-1}, i\in\{1,\dots,n\}$ we denote by $C(g)=
B_i, \hbox { and } C'(g)=A_i$ if $g=f_i$ and $C(g)= A_i\hbox{ and }
C'(g)=B_i$ if $g=f_i^{-1}$, so that for every $g$ one has
$g(\cE\setminus C'(g))\subset C(g)$.
Note that if $g_1\neq g_2^{-1}$ then $g_2(C(g_1))\subset C(g_2)$ so that $g_2\circ g_1(\cE\setminus C'(g_1))\subset C(g_2)$.

\begin{lemm} Let $(\cE,d,\{f_i\})$ be a compact contracting ping-pong.
For every $\varepsilon>0$ there is $\ell\in \NN$ such that for every
reduced word $g_\ell\circ\cdots \circ g_1$,
$g_i\in\{f_i, f_i^{-1}, i\in\{1,\dots,n\}\}$ one has

$$diam( g_\ell\circ \cdots\circ g_1(\cE\setminus C'(g_1))) <\varepsilon$$
\label{l.contract}
\end{lemm}
\begin{demo}
Using the compacity of the set of points $x,y$ such that $d(x,y)\geq\varepsilon$ we get that 
there is $0<\delta<1$ such that if $x,y\in \cE\setminus C'(g)$, and $d(x,y)\geq \varepsilon$ then $d(g(x),g(y))\leq \delta\cdot d(x,y)$. 
\end{demo}

Let $\Si_0=\{f_i, f_i^{-1}, i\in\{1,\dots,n\}\}^{Z}$ be the set of infinite words with letters equal to $f_i^{\pm 1}$, endowed with the product topology. An infinite word $(g_i)_{i\in\ZZ}$ is called {\em reduced} if for any $n$ the finite word $(g_i)_{-n<i<n}$ is reduced.
We denote by $\Si=\{(g_i)\in\Si_0, (g_i) \mbox{ is reduced} \}$ the 
subspace of reduced words,
 $\De= \Si\times \cE$ and $\Pi\colon \De\to \Si$ the natural projection. 
Denote by $\sigma$ the shift on $\Si$, that is $\sigma(g_i)= (h_i)$ where $h_i= g_{i+1}$, and by $\tilde \sigma$ the map on $\De$ defined by  $\tilde\sigma ((g_i),x))= (\sigma(g_i),g_0(x))$. One verifies easely that $\sigma$ and $\tilde\sigma$ are homeomorphisms. Notice that $\tilde\sigma$ is a multiplicative cocycle over $\sigma$.

\vskip 2mm
The topological picture of the ping-pong may be
completely understood:

\begin{prop}
With the notation above, there are exactly two continuous sections $s^+\colon \Si\to \De$ and $ s^-\colon\Si\to \De$ which are $\tilde\sigma-$invariant. Moreover, $s^+(\Si)$ is a topological attractor for $\tilde \sigma$  whose basin is $\De- s^-(\Si)$ and $s^-(\Si)$ is a topological repellor for $\tilde \sigma$ with basin $\De- s^+(\Si)$ and these two sections are disjoint.
\label{4.sections}
\end{prop}
\begin{demo}
Let $(g_i)\in\Si$ be a reduced word. For every $n\in\NN$, consider the compact sets 
$$K_n^+ = g_{-1}\circ g_{-2}\circ\cdots\circ g_{-n}(\cE\setminus C'(g_{-n})\subset C(g_{-1}) $$
and
$$K_n^- = g_{0}^{-1}\circ g_1^{-1}\circ\cdots\circ g_{n-1}^{-1}(\cE\setminus C(g_{n-1})\subset C'(g_0) $$
Using the fact that the word $(g_i)$ is reduced, one shows easily that these sequences of compact sets are decreasing with $n$: $K_{n+1}^+\subset K_n^+$ and $K_{n+1}^-\subset K_n^-$.
Moreover as $g_0\neq g_{-1}^{-1}$ one has $C(g_{-1})\cap C'(g_0)=\emptyset$, so that $K_n^+\cap K_n^-=\emptyset$. Finally,  Lemma~\ref{l.contract} ensures that the diameter of $K_n^+$ and $K_n^-$ goes uniformly to $0$. 
We define then
$$
s^-((g_i))=\bigcap_{n\in\NN} K_n^- \quad\mbox{and}\quad s^+((g_i))=\bigcap_{n\in\NN} K_n^+
$$
\end{demo}

\subsection{Schottky Groups}
A {\em Schottky group} of rank $n$ is a finitely generated group $\Ga\subset PSL(2,\CC)$ having $2n$ disjoint circles $C_1, C'_1,\dots, C_n, C'_n$ bounding a domain $D\subset \CC P^1=\CC\cup\{\infty\}$, and a system $g_1,\dots,g_n$ of generators such that $g_i(C'_i)=C_i$ and $g_i(D)\cap D=\emptyset$ (see \cite{Ma}). Using the discs $A_i,B_i$ bounded by the circles $C_i,C'_i$ respectively and disjoint from $D$, one see that $\Ga$ is a ping-pong group of $Aut(\CC P^1)$, moreover it is a compact contracting ping-pong group.

\subsection{Geodesics and Reduced Words}

\begin{lemm} Let $S$ be a finite non-compact hyperbolic Riemann
surface,
endowed with its natural hyperbolic metric.
There are $\gamma_1,\dots,\gamma_k$ complete mutually disjoint geodesics
whose ends arrive to punctures of $S$, such that
the complement $S\setminus\bigcup_1^k \gamma_i$
is connected and simply connected, the $\gamma_i$ bound a
fundamental domain of $S^\prime$ in its universal cover $\DD$ and
 the fundamental domain is a $2k$ sided
polygon whose vertices are on the circle at infinity of $\DD$.
\end{lemm}
\begin{demo}
Let $\beta_1,\dots,\beta_k$ be a maximal set of
non-homotopic mutually disjoint curves
whose ends arrive to punctures of $S$.
Clearly, by removing them from $S$ we obtain a connected
simply connected domain (for otherwise we could pick
and additional $\beta_{k+1})$.
Lift them to the universal cover of $S$ and replace the lifts
of $\beta_j$ by the geodesics that have the same endpoints.
Pushing down these geodesics to $S$, gives the desired curves $\gamma_i$.
\end{demo}

Now  fix an origin $x_0\in S\setminus \bigcup_1^k\gamma_i$.
For each $i$ there is a unique geodesic segment $\alpha_i$ joining $x_0$ to $x_0$ and cutting $\gamma_i$ at exactly one point, with the positive orientation, and not
cutting $\gamma_j,\quad j\neq i$.

\begin{lemm} The closed paths $\alpha_i$ build a system of generators of the fundamental group
$\pi_1(S,x_0)$. More precisely the fundamental group is the free group generated by the $\alpha_i$.
\end{lemm} 
\begin{demo} The union of the $\alpha_i$ is a bouquet
 of circles and we verify easily that $S$ admits a retraction by
 deformation on this bouquet of circles.
\end{demo}

Now  fix an orientation on each geodesic $\gamma_i$ and  call
 $\gamma_i$ the oriented geodesic. Given any vector $u\in T^1_xS$ at a point $x\in S\setminus\bigcup_1^k\gamma_i$,
the geodesic $\gamma_u$ has two possibility:
\begin{enumerate}
\item either one of its ends goes to one puncture of $S$,
\item or $\gamma_u$ cuts transversely infinitely many times (in the future and in the past) the geodesics $\gamma_i$. 
\end{enumerate}

\begin{defi} {\em The itinerary of the geodesic $\gamma_u$} is
the sequence $b(u)=(b_i)_{i\in\ZZ}$ defined as follows: 

\vskip 2mm
$b_i$ is $\alpha_i^{\pm 1},\quad i\in\{1,\dots,k\}$ if the $(i-1)^{th}$ intersection of $\gamma_u$ with $\bigcup \gamma_l$ belongs to $\gamma_i$ and the coefficient is $+1$ or $-1$ according if the orientation of $\gamma_u$ followed by the orientation of $\gamma_i$ is a direct or inverse basis of the tangent space.  
\end{defi} 

\begin{lemm} For any $u
\in T^1_xS^\prime$ the itinerary $b(u)$ is a (finite or infinite)
reduced word in the letters $\alpha_i^{\pm 1}$,
where $b_0$ corresponds to the first intersection point.
\end{lemm}
\begin{demo} If a  segment in the fundamental domain cuts 2
times the same geodesic $\gamma_i$ with opposite direction,
then its lift on $\DD$ will cut 2 times the same lift of $\gamma_i$. So this segment cannot be geodesic.
\end{demo}

Given the geodesic $\gamma_u$, and a time $t_0\in\RR$
such that $\gamma_u([0,t_0])\notin\bigcup_1^k \gamma_i$, we get a closed path $\tilde\gamma_u(t)$ joining respectively $\gamma_u(0)$ and $\gamma_u(t)$ by a geodesic segment in the fundamental domain. Moreover if $t>0$ and if the segment $\gamma_u([0,t])$ cuts $\ell+1$ times the geodesic $\gamma_i$,  then the
closed path $\tilde \gamma_u(t)$ is homotopic to $\beta_0\cdot \beta_1\cdots\beta_\ell$ where $\beta_j$ is a closed path $\alpha_i^{\pm 1}$ according to the letter $b_j= \alpha_i^{\pm 1}$.

\begin{coro}
The geodesic $\gamma_u$ defines a (finite or infinite) reduced word in $\pi_1(S,x_0)$ for the basis $\alpha_i,\quad i\in\{1,\dots,k\}$.
\end{coro}

\subsection{Proof of  Theorem~\ref{t.Schottky}}

Let $G_e$ and $G_f$ be free groups generated by
$e=\{e_1,\dots,e_k\}$
and $f=\{f_1,\dots,f_\ell\}$, respectively.
Denote by $\Ga_e$ and $\Ga_f$ their Cayley graphs for the
given basis. Both Cayley graphs are trees.
Let $\rho:G_e \rightarrow G_f$ be a group isomorphism.
Any infinite word $b=(b_j)_{j\in \ZZ},\quad b_j\in\{e_i^{\pm1}, 1\leq i\leq k\}$
defines an infinite path $\sigma(b)$ in the Cayley graph $\Ga_e$. This path $\sigma(b)$ is a geodesic if and only if the word $b$ is reduced
(see [12] for background material on hyperbolic groups).

\begin{defi} We say that an infinite path $\sigma\subset \Ga_e$ is 
stretchable
if it is properly embedded (namely,
only a bounded part of the path remains in a given compact
set of the Cayley graph). It is strictly stretchable if its 2
ends correspond to two distinct ends $\sigma_-$ and $\sigma_+$ of the Cayley graph. The unique geodesic joining $\sigma_-$ to $\sigma_+$ is the reduction $\sigma^r$ of $\sigma$.
\end{defi}

\vskip 2mm
 \begin{lemm} Let $b$ an infinite word in the letters $(e_i)$. Let 
$c:=\rho(b)$
be the corresponding word in the letters $f_i$. Then $b$ is stretchable if and only if $c$ is stretchable.
$\rho$ induces a homeomorphism from the
boundary of $\Gamma_e$ to the boundary of $\Gamma_f$
by associating to the boundary point $b$ the
boundary point $\rho(b)$.
\end{lemm}

\begin{demo} 
Given any word $b$ in the letters $e_i$, $\rho$ produces
 a reduced word $c:=\rho(b)$ in the
letter $f_i$ obtained as follows:
Change each letters $b_j=e_i^{\pm1}$ by
the reduced word $\rho(b_j)$ written in terms of $f$.
Do the appropiate cancellations to obtain the reduced word
$c$. By [12] p.7, the isomorphism $\rho$ induces a
quasi-isometry of the Cayley graphs, hence $b$ is stretchable
if and only if $c$ is. 
\end{demo}

A stretchable word $a$ in a free group defines two points $a_-$ and $a_+$ in the boundary of the group. So there is a unique geodesic $c^r$ in the Cayley graph of the group, 
which corresponds to a reduced word on the group, joining  $a_+$  to $a_-$.  
Using the same notation as in Proposition 6.4,
define $s^-(a) =s^-(c^r)$ and $s^+(a)=s^+(c^r)$.
Denote by $\hat\Si$ the set of stretchable 
infinite words whose letters are
the generators of the Schottky group $G$. The reduced word corresponding to $c^r$ above 
belongs to $\hat\Si$. 
Recalling that in this case, the group acts on $\CC P^1$;  $\sigma$ is the shift 
on 
$\hat \Si$, being a homeomorphism, because $\hat \Si$  is $\sigma$ invariant. Recall that $\tilde \sigma$ is the map on
$\hat\Sigma\times \CC P^1$ defined by $\tilde\sigma(a,x)=(\sigma(a),a_0(x))$.
Then, since the Schottky group defines a compact contracting ping pong, Proposition 6.4
implies immediately the following:

\begin{lemm} Let $a$ be a stretchable word  in a Schottky group $G\subset SL(2,\CC)$ and $b$ the image of $a$ by the shift. Then $s^\pm(b)=a_0(s^\pm(a))$.
 The map $s^\pm\colon a\mapsto (a,s^\pm(a))$ defines an $\tilde\sigma-$ measurable section of the trivial fibration $\hat\Si\times\CC P^1\to \hat\Si$.
\end{lemm}
\qed

\vskip 2mm
\noindent
{\bf Proof of Theorem~\ref{t.Schottky}:}
Let 
$\rho\colon
\pi_1(S,x_0)\to SL(2,\CC)$
be an injective representation 
with $G=\rho(
\pi_1(S,x_0))$
 a Schottky group.
Notice that the set of vector $u\in T^1S$  such that the
corresponding geodesic $\gamma_u$ goes to a puncture of $S$
has zero Lebesgue measure.
 
For any  unit vector $u$ at a point of the  fundamental domain
such that the geodesic $\gamma_u$ has no end at a puncture of $S$,
the word $\rho(b(u))$ is a stretchable word of the Schottky group.
For any point $x$ of the fundamental domain of  $S$ we denote by $H_x$
the holonomy of the foliation $\cF_\rho$ from the fiber
over $x$ to the fiber over $x_0$ by a path contained inside the
fundamental domain. This holonomy is well defined because the
fundamental domain is simply connected. 
So we define $s^\pm\colon T^1S\to T^1\cF_\rho$
as  $s^\pm(u)=H_x^{-1}(s^\pm(\rho b(\gamma_u)))$.
By construction the sections $s^\pm$ are defined
Liouville almost everywhere,
are measurable, and are the sections of largest expansion
and contraction.
The continuity of $s^\pm$ follows from
the topological way of constructing the sections
in Proposition 6.4 and the fact that the map which associates
the point at infinity of the Cayley graph of the presentation
of $\pi_1(S)$ to the point at infinity of the
Cayley graph of the Schottky group is continuous,
by Lemma 6.11.
This proves Theorem~\ref{t.Schottky}.
\qed

\vskip 2mm

\noindent
{\bf Remark:}
Observe that Schottky representations over punctured Riemann surfaces
never satisfy the integrability condition $(3)$ due to Theorem 3, since
all its elements are hyperbolic and so, in particular, the maps
corresponding to loops around a puncture. By the way we chose
the presentation of the fundamental group (Lemma 6.5) the geodesics
give rise to reduced words. Assume now that the image under
$\rho$ of these generators of $\pi_1(S)$ are generators
of the Schottky group, then we will have that
there are no cancellations in the words corresponding to $\rho$(geodesic).
For  the general geodesic in $S$,
the ratio between the number of letters to the length of the
geodesic goes to infinity as the length of the geodesic
goes to infinity, since by ergodicity of
the geodesic flow the average time that the general geodesic
spends in a small disk around the puncture is proportional
to the area of the disk and the number of turns that the geodesic
does around the pucture is $cot(\eta)$ by Proposition 5.4.
This shows that the `Lyapunov exponents' of these Schottky representations
are $\pm\infty$.

\vskip 2mm
\noindent
{\bf Remark:} 
If $S$ is compact and the group $\tilde\rho(\pi_1(S))$ 
is non-cyclic but contained in  a Schottky
group, it follows from the results in [3] 
that there are 
positive and negative Lyapunov exponents, and hence
sections of largest expansion and contraction, 
but they will only be measurable sections now
due to cancellations in the reduced words.

  \section{Foliation Associated to the Canonical  Representation}

\subsection{The Geometry of the Bundles}

Let $S$ be a hyperbolic Riemann surface, and denote by $\pi\colon \HH^+\to S$ its universal cover
by the upper half plane $\HH^+$. Fix a point  $x_0\in S$, and $\bar x_0\in\pi^{-1}(x_0)$. Denote by 
$$\rho_{can}\colon\pi_1(S,x_0)\to PSL(2,\RR)\subset PSL(2,\CC)$$ the covariant representation obtained 
by the covering transformations.   We consider now the suspension foliation $\cF_{can}$ associated to the representation 
$\rho_{can}$ (that is a foliation in $M_{can}$ whose holonomy is given by $Hol(\gamma)=\rho_{can}(\gamma)^{-1}$).  
 
 \begin{defi} 
The representation $\rho$, the $\CC P^1$
  bundle $M_{can}$ and the foliation $\cF_{can}$  are called the {\em  canonical
 representation, bundle and foliation of the hyperbolic Riemann surface $S$}.  
\label{d.tautologic}
\end{defi}

Denote by $\iota\colon\HH^+\to\CC P^1$ the usual inclusion
of the upper half plane in the projective line. 
We have the canonical action
$$\pi_1(S,x_0) \times [\HH^+ \times \CC P^1]
\longrightarrow
[\HH^+ \times \CC P^1]
\hskip 1cm, \hskip 1cm
(\gamma,x,z) \rightarrow (\rho_{can}(\gamma)(x),
\rho_{can}(\gamma)(z))
$$
corresponding to the representation
$$\rho_{can}\times\rho_{can}: \pi_1(S,x_0)
 \rightarrow PSL(2,\RR)\times
PSL(2,\CC)$$
The quotient $\Pi:M_{can} \rightarrow \HH^+/\rho_{can} = S$
is a 2-dimensional complex manifold and the projection
to the first factor gives it the structure of a $\CC P^1$
bundle over $S$.

\vskip 2mm
For any $\alpha\in PSL(2,\RR)$ one has $\iota\circ \alpha_{\HH^+}=\alpha_{\CC P^1}\circ\iota$. 
Denote by $\tilde \De$ the diagonal $\tilde\De=\{(z,\iota(z))| z\in\HH^+\}$. Then for each $\gamma\in\pi_1(S,x_0)$ and each $z\in\HH^+$ one gets:
$$ (\rho_{can}(\gamma)z,\rho_{can}(\gamma)\iota(z))= (\rho_{can}(\gamma)z, \iota(\rho_{can}(\gamma)(z))\in\tilde\De,$$ 
so the diagonal $\tilde\De$ is  invariant
 by the action   of 
$\rho_{can}\times\rho_{can}$
and induces  in the complex surface $M_{can}$ a Riemann surface $\De$ and the projection $\Pi$ induces a biholomorphism  $\De\to S$. The diagonal
$\De$ is the image of a holomorphic section of the bundle $M_{can} \rightarrow S$. 

As the representation $\rho_{can}$ has its values in $PSL(2,\RR)$, the circle bundle $\HH^+\times \RR P^1$ is invariant by the action of $\rho(\gamma),\quad \gamma\in\pi_1(S,x_0)$, so that it 
defines  $M_{can}^\RR \subset M_{can}$ an $\RR P^1-$subbundle.
For every point $p$ of $S$ we will denote by
 $\RR P^1_p\subset \CC P^1_p$ the fiber of these bundles over $p$.
$M_{can}^\RR$  is disjoint from the diagonal $\De$. 

Consider now the unit tangent spaces $\Pi_*\colon T^1\cF_\rho\to T^1S$.
Notice that every unit vector $u$ at a point $p\in S$ lifts canonically to a unit vector tangent to $\cF$ at any point $\tilde p$ in the fiber $\CC P^1_p$. So the diagonal $\De$ induces canonically a section 
$\De_* \colon T^1S\to T^1\cF$:
$$
\begin{matrix}
M_{can} &
\leftarrow &
T^1\cF \cr
\Delta\uparrow\downarrow \Pi && \Pi_*\downarrow\uparrow 
 \Delta_*\cr S &
\leftarrow &
T^1S \cr
\end{matrix}
$$

\begin{defi} For every unit
 vector $u \in T^1_p\HH^+$, the geodesic $\gamma_u$ through $p$ tangent to $u$ has 
its extremities $\tilde \sigma^+(u)$ and $\tilde\sigma^-(u)$ in $\RR P^1$. 
This defines 2 smooth sections $\tilde\sigma^\pm:T^1\HH^+ 
\rightarrow T^1\HH^+ \times \CC P^1$.
Let $Y_u$ be the 
 holomorphic vector field  
on $\CC P^1$ vanishing at
$\tilde\sigma^\pm(u)$
and having $Y_u(p) = u$.
Let $\tilde Y$ be the smooth vector field
defined on $T^1\HH^+ \times \CC P^1$ by $\tilde Y(v,.) := Y_v(.)$.
$\tilde Y$  is tangent to 
the fibers $\{u\}\times\CC P^1, \quad u\in T^1\HH^+$.
\label{r.y}
\end{defi} 

 Note that if $\tilde\sigma^-(u)= 0 \in \CC P^1$, $\tilde\sigma^+(u)=\infty$ 
and $u$ is the vector $i \in T_i\HH^+$  then 
$Y_u$ is the vector field $z\frac{\partial}{\partial z}$. So for every $u$, $Y_u$ is conjugate to $z\frac{\partial}{\partial z}$.
The hyperbolic norm of $Y_u$ along
the geodesic $\gamma_u$ is uniformly $1$. So the flow of $Y_u$ induces the translations along this geodesic.
The derivative of $Y_u$ at the point $\tilde\sigma_0^-(u)$ is equal to $1$, and this does not depend on the metrics on $\CC P^1$.
The flow lines of the vector field
  $z\frac{\partial}{\partial z}$
consist of semirays through $0$ having a north to south
pole dynamics, with $0$ as a hyperbolic repellor and $\infty$ as a 
hyperbolic attractor. The vertical ray is a geodesic in $\HH^+$.

 \begin{lemm} The sections $\tilde\sigma^+$ and $\tilde \sigma^-$ 
 and the vector field $\tilde Y$ 
are invariant by every $
T \in PSL(2,\RR)$, i.e.:
$$\sigma^\pm(T_*(v)) = T(\sigma^\pm(v)) \hskip 1cm,
\hskip 1cm (T_*\times T)_*\tilde Y = \tilde Y$$
 \label{l.sections} 
\end{lemm} 
\begin{demo}
The endpoints of the geodesic determined by $T_*v$
are $T(\sigma^\pm(v))$, so they are invariant,
as well as $Y_{T_*(v)} = T_*Y_v$, by its definition.
\end{demo}

 The sections $\tilde\sigma^\pm$   induce in the quotient bundle sections $\sigma^\pm$
from $T^1S$ to the $\RR P^1-$subbundle of $T^1\cF$,
and $\tilde Y$ induces a vector field $Y$ on $T^1\cF$.
The sets $\sigma^\pm(T^1S)$ are the zero sets of $Y$.
  
\begin{coro} The diagonal $\De$, $\sigma^+$ and $\sigma^-$ are $3$ smooth sections of $T^1\cF \rightarrow T^1S$, 
pairwise disjoint, and hence define a smooth trivialisation of the $\CC P^1-$fiberbundle 
$$[T^1\cF\to T^1S ]\sim [T^1S\times \CC P^1\to T^1S]$$ 
sending $\sigma^+$ to $\infty$, $\sigma^-$ to $0$ and $\De$ to $1$.
\label{c.trivial}
\end{coro} 
\begin{demo} The unique thing we need to prove is that the sections are two by two disjoint.
 $\sigma^+$ and $\sigma^-$ are included in the $\RR P^1$ bundle which is disjoint from $\De$, since the image of 
$\Delta$ is in the upper half plane.
 The 2 points $\sigma^\pm(u)$ 
are the extremities in $\RR P^1$ of 
 a geodesic in $\HH^+\subset \CC P^1$, so they are different.
\end{demo}

We will denote by $|\cdot|$ the Fubini Study metric on the fibers of $T^1\cF_{can}$ induced by the trivialisation $T^1\cF=T^1S\times \CC P^1$ given by Corollary~\ref{c.trivial}.

\begin{rema}In the trivialisation $T^1\cF_{can} \sim T^1S\times \CC P^1 $ given by Corolary~\ref{c.trivial} the flow $Y$ admits the sections $T^1S\times\{0\}$ and $T^1S\times\{\infty\}$ as zeros and the vertical derivative on every point $(u,0)$ is $1$.
So in this coordinates the vector field $Y$ is $(0,z\frac{\partial}{\partial z})$. 
\end{rema}

\subsection{The Foliated Geodesic Flow} 
Denote by $X$ and $X_{can}$ the
infinitesimal generators of the geodesic and the
foliated geodesic flows on $T^1S$ and $T^1\cF_{{can}}$, respectively,
and $\varphi$ and $\Phi$ the corresponding flows, as in $(2.2)$.

\begin{prop}
The vector fields $X_{can}$ and $Y$ on $T^1\cF_{{can}}$  commute.
In particular, the set $Zero(Y)$ is invariant by $X$, so that $\sigma^+$ and $\sigma^-$ are invariant by $X$.
\label{p.commute}
\end{prop}
\begin{demo}
It suffices to show that $\Phi_{t*}Y=Y$, since
$$[X,Y] = \lim_{t \rightarrow \infty} \frac{1 }{ t}
[\Phi_{t*}Y - Y] = 0.$$
The proof of this is easier on the universal cover $T^1\HH^+\times \CC P^1$. 
Let $\tilde X$ be the lift of $X$ to 
the universal covering space $T^1\HH^+ \times \CC P^1$. 
In this trivialisation, the foliated geodesic flow is generated by 
$(X,0)$, 
 since the foliation is horizontal. 
So it is enough to prove the following statment:
\begin{claim}
Let $u$ and $v$ be  unit vectors tangent to the same geodesic $\gamma$ of $\HH^+$ at  $x$ and $y$, and inducing the same orientation of $\gamma$. Then the vector fields $Y_u$ and $Y_v$ on $\CC P^1$ coincide. 
\end{claim}
To prove the claim it is enough to notice that $\iota_*(u)$ and $\iota_*(v)$ are unit vectors
for the hyperbolic metric of $\HH^+\subset \CC P^1$ tangent at the points $\iota(x)$ and $\iota(y)$ to the geodesic (for the hyperbolic metric) $\iota(\gamma)$. 
The vector field $Y_u$ is tangent to every point of
$\gamma_u$ and its hyperbolic norm is $1$, moreover the orientation induced by $Y_u$ on $\gamma$ cannot change. So $Y_u(y)=\iota_*(v)$ and so $Y_u=Y_v$. Hence $\Phi_{t*}Y=Y$ as required.

The claim shows that for every $u$ and every $v=\phi_t(u)$
  the vertical vector field $Y$ on $\{v\}\times\CC P^1$ is $\Phi_{t*}(Y|_{\{u\}\times \CC P^1})$.
  Hence $\Phi_{t*}Y=Y$ as required.
\end{demo}

\begin{prop} The vector field $Z=X+Y$ is tangent to the diagonal $\tilde\De$.
\label{p.tangent}
\end{prop}
\begin{demo} The proof is easier on the  cover $\HH^+\times \CC P^1$.
Consider the following diagram:
$$
\begin{array}{ccc}
                               
         &\tilde p&\\           
\HH^+\times \CC P^1 
 & 
\leftarrow 
& 
T^1 {\tilde\cF}=T^1\HH^+\times\CC P^1
\\
 \Pi \downarrow\quad\uparrow \De
&
&
 \Pi_*\downarrow\uparrow \De_* 
\\
{\HH^+} 
& 
\leftarrow 
& 
T^1\HH^+ 
\\
		&p&\\
\end{array} 
$$
 To show that $X+Y$ is tangent to the diagonal $\tilde \De$ it is enough to show that, for every $u_x\in T^1\HH^+, x\in\HH^+$ the vector $\tilde p_*((X+Y)(u_x,\iota(x))$ is tangent to $\De$ at the point $(x,\iota(x))$.
On one hand, $\tilde p_*(X(u_x,y))$ is the horizontal vector $(u_x,0)$ at the point $(x,y)$. 
On the other , $\tilde p_*(Y(u_x,\iota(x)))$
 is the vertical vector $(0,\iota_*(u_x))$ at the point $(x,\iota(x))$. 
So the vector $\tilde p_*((X+Y)(u_x,\iota(x)))$ is the vector $(u_x,\iota_*(u_x))$ at the point $(x,\iota(x))$ and is tangent to $\De$.
\end{demo}

\begin{coro}The flow $Z_t$ of $Z$ is horizontal in the trivialisation $T^1\cF$. In particular it induces isometries on the fibers $\CC P^1$ endowed with the metric $|\cdot|$.
\end{coro}
\begin{demo} As $X$ and $Y$ commute and all preserve the fibration so does $Z$. Moreover, as $X$ and $Y$ induce on the fiber maps belonging to $SL(2,\RR)$ so does $Z$.
 To prove the corollary its suffices to show
 that $Z$ preserves the $3$ sections 
$\tilde\De$, $\tilde\sigma^+$ and $\tilde\sigma^-$. 
$Z$ is tangent to $\tilde \De$ by Proposition 7.7.
$Y$ vanishes on $\sigma^\pm(T^1S)$ and $X$ is tangent to them
by Proposition 7.6.
\end{demo}

\noindent {\bf Proof of Theorem~\ref{t.tautologic}}
The foliated geodesic flow is $X=Z-Y$. As these 
flows commute $X_t= Y_{-t}\circ Z_t$,
where the notation corresponds to
the flows of the corresponding
vector fields. In the trivialisation given by Corollary~\ref{c.trivial} the flow $Z_t$ induces the identity on the fibers and $Y_{-t}$ is the homothety $z\to e^{-t}z$.
Hence we obtain a contraction in the projective space,
which may be translated to the affine space. This
means that there is a section of largest expansion
and contraction. The sections are smooth sections.
The geodesic flow is recurrent hence the $\omega$ limit set of any point
not in $\sigma^-(T^1S)$ is contained in $\sigma^+(T^1S)$.
The $\alpha$ limit set of any point
not in $\sigma^+(T^1S)$ is contained in $\sigma^-(T^1S)$.
Along $\sigma^\pm(T^1S)$ the foliated geodesic flow $X_\De$
is hyperbolic. This proves the Theorem~\ref{t.tautologic}.
\qed

\subsection{ Representation Topologically Equivalent to the Canonical Representation}

Let $$V := \{\rho := (A_1,\ldots,A_g) \in PSL(2,C) \ / \
\Pi_{1}^g[A_{2i-1},A_{2i}] = Id \}$$
be the complex algebraic variety parametrizing 
representations of the fundamental group
$\pi_1(S)$ of the compact Riemann surface of genus $g \geq 2$,
where $[A,B]:=ABA^{-1}B^{-1}$.
We also have an action 
$$PSL(2,C) \times V \rightarrow V
\eqno (7.1)$$
given by conjugation.
Let $\rho_0$ be the representation corresponding to the
canonical representation.  Bers's simultaneous 
uniformisation ([19]) implies
that there is an open connected
set $U \subset V$ containing $\rho_0$ such that all representations
in $U$ are quasiconformally conjugate, and 
there is a surjective 
map $$U \rightarrow Teich^g \times Teich^g$$
which associates to each representation $\rho \in U$ the
Riemann surfaces obtained by quotienting the region of discontinuity of
$\rho$ by $\rho$, and its
fibers are the $PSL(2,\CC)$ orbits $(7.1)$.

 \begin{prop} For any representation $\rho$ in the above open set $U$,
 the Riccati equation with monodromy $\rho$
 has a unique SRB-measure with basin of attraction of total
 Lebesque measure
for positive and for negative times.
 \end{prop}
\begin{demo}
By Theorem 6, the assertion is true for the canonical representation
$\rho_{can}$.
By Bers's simultaneous uniformization, there is a
quasiconformal map $h:\CC P^1 \rightarrow \CC P^1$ conjugating
the action of
$\rho_{can}$  to the action of $\rho \in U$.
We may use this map to obtain
a homeomorphism over $T^1S$ of the $\CC P^1$-bundles
$H:Proj(E_{{can}}) \rightarrow Proj(E_\rho)$
conjugating the geodesic flows. This homeomorphism
is absolutely continuos, since horizontally it is the
identity and vertically it is the quasiconformal map $h$,
which is absolutely continuos. Hence
$ Proj(E_\rho)$ has a unique SRB-measure for positive and negative times,
and it is $H_*(\mu^\pm)$.
\end{demo}

\vskip 2mm
Christian Bonatti ( bonatti@@u-bourgogne.fr)

Laboratoire de Topologie, UMR 5584 du CNRS

B.P. 47 870,
21078 Dijon Cedex, France

\vskip 2mm
Xavier G\'omez-Mont ( gmont@@cimat.mx)

Ricardo Vila ( vila@@cimat.mx)

CIMAT

A.P. 402,
Guanajuato, 36000, M\'exico

\end{document}